\documentclass[journal]{IEEEtran}

\usepackage[mathscr]{eucal}
\usepackage[cmex10]{amsmath}
\usepackage{epsfig,epsf,psfrag}
\usepackage{amssymb,amsmath,amsthm,amsfonts,latexsym}
\usepackage{amsmath,graphicx,bm,xcolor,url}
\usepackage[caption=false]{subfig} 
\usepackage{fixltx2e}
\usepackage{array}
\usepackage{verbatim}
\usepackage{bm}
\usepackage{algorithmic, cite}
\usepackage{algorithm}
\usepackage{verbatim}
\usepackage{textcomp}
\usepackage{mathrsfs}
\usepackage{epstopdf}

\catcode`~=11 \def\UrlSpecials{\do\~{\kern -.15em\lower .7ex\hbox{~}\kern .04em}} \catcode`~=13 

\allowdisplaybreaks[3]

\newcommand{\calI}{\mathcal{I}}

\newcommand{\calN}{\mathcal{N}}
\newcommand{\calO}{\mathcal{O}}

\newcommand{\calU}{\mathcal{U}}

\newcommand{\bA}{\mathbf{A}}
\newcommand{\bb}{\mathbf{b}}
\newcommand{\bB}{\mathbf{B}}

\newcommand{\bC}{\mathbf{C}}
\newcommand{\bd}{\mathbf{d}}
\newcommand{\bD}{\mathbf{D}}
\newcommand{\be}{\mathbf{e}}
\newcommand{\bE}{\mathbf{E}}

\newcommand{\bG}{\mathbf{G}}

\newcommand{\bH}{\mathbf{H}}

\newcommand{\bL}{\mathbf{L}}

\newcommand{\bN}{\mathbf{N}}

\newcommand{\bO}{\mathbf{O}}

\newcommand{\bP}{\mathbf{P}}
\newcommand{\bq}{\mathbf{q}}
\newcommand{\bQ}{\mathbf{Q}}

\newcommand{\bR}{\mathbf{R}}
\newcommand{\bs}{\mathbf{s}}
\newcommand{\bS}{\mathbf{S}}

\newcommand{\bU}{\mathbf{U}}

\newcommand{\bV}{\mathbf{V}}

\newcommand{\bz}{\mathbf{z}}

\DeclareMathAlphabet{\mathbsf}{OT1}{cmss}{bx}{n}
\DeclareMathAlphabet{\mathssf}{OT1}{cmss}{m}{sl}

\DeclareSymbolFont{bsfletters}{OT1}{cmss}{bx}{n}  
\DeclareSymbolFont{ssfletters}{OT1}{cmss}{m}{n}
\DeclareMathSymbol{\bsfGamma}{0}{bsfletters}{'000}
\DeclareMathSymbol{\ssfGamma}{0}{ssfletters}{'000}
\DeclareMathSymbol{\bsfDelta}{0}{bsfletters}{'001}
\DeclareMathSymbol{\ssfDelta}{0}{ssfletters}{'001}
\DeclareMathSymbol{\bsfTheta}{0}{bsfletters}{'002}
\DeclareMathSymbol{\ssfTheta}{0}{ssfletters}{'002}
\DeclareMathSymbol{\bsfLambda}{0}{bsfletters}{'003}
\DeclareMathSymbol{\ssfLambda}{0}{ssfletters}{'003}
\DeclareMathSymbol{\bsfXi}{0}{bsfletters}{'004}
\DeclareMathSymbol{\ssfXi}{0}{ssfletters}{'004}
\DeclareMathSymbol{\bsfPi}{0}{bsfletters}{'005}
\DeclareMathSymbol{\ssfPi}{0}{ssfletters}{'005}
\DeclareMathSymbol{\bsfSigma}{0}{bsfletters}{'006}
\DeclareMathSymbol{\ssfSigma}{0}{ssfletters}{'006}
\DeclareMathSymbol{\bsfUpsilon}{0}{bsfletters}{'007}
\DeclareMathSymbol{\ssfUpsilon}{0}{ssfletters}{'007}
\DeclareMathSymbol{\bsfPhi}{0}{bsfletters}{'010}
\DeclareMathSymbol{\ssfPhi}{0}{ssfletters}{'010}
\DeclareMathSymbol{\bsfPsi}{0}{bsfletters}{'011}
\DeclareMathSymbol{\ssfPsi}{0}{ssfletters}{'011}
\DeclareMathSymbol{\bsfOmega}{0}{bsfletters}{'012}
\DeclareMathSymbol{\ssfOmega}{0}{ssfletters}{'012}

\DeclareMathOperator{\spn}{span}

\newtheorem{theorem}{Theorem} 
\newtheorem{lemma}[theorem]{Lemma}

\newtheorem{remark}{Remark}

\newtheorem{assumption}{Assumption}

\newcommand{\qednew}{\nobreak \ifvmode \relax \else
      \ifdim\lastskip<1.5em \hskip-\lastskip
      \hskip1.5em plus0em minus0.5em \fi \nobreak
      \vrule height0.75em width0.5em depth0.25em\fi}

\begin{document}
\title{ { High Dimensional Low Rank plus Sparse Matrix Decomposition}}

\author{Mostafa~Rahmani, \IEEEmembership{Student Member,~IEEE} and George~K.~Atia,~\IEEEmembership{Member,~IEEE} 
\thanks{This material is based upon work supported by NSF CAREER Award CCF-1552497 and NSF grant No. CCF-1320547.

The authors are with the Department of Electrical and Computer Engineering, University of Central Florida, Orlando, FL 32816 USA (e-mail: mostafa@knights.ucf.edu, george.atia@ucf.edu).}
}

\markboth{}%
{Shell \MakeLowercase{\textit{et al.}}: Bare Demo of IEEEtran.cls for Journals}
\maketitle

\begin{abstract}
This paper is concerned with the problem of low-rank plus sparse matrix decomposition for big data.
Conventional algorithms for matrix decomposition use the entire data to extract the low-rank and sparse components, and are based on optimization problems with complexity that scales with the dimension of the data, which limits their scalability. 
Furthermore, existing randomized approaches mostly rely on uniform random sampling, which is quite inefficient for many real world data matrices that exhibit additional structures (e.g. clustering). In this paper, a scalable subspace-pursuit approach that transforms the decomposition problem to a subspace learning problem is proposed. The decomposition is carried out using a small data sketch formed from sampled columns/rows. Even when the data is sampled uniformly at random, it is shown that the sufficient number of sampled columns/rows is roughly $\calO(r \mu)$, where $\mu$  is the coherency parameter and $r$ the rank of the low-rank component. In addition, adaptive sampling algorithms are proposed to address the problem of column/row sampling from structured data.
 \textcolor{black}{We provide an analysis of the proposed method with adaptive sampling and show that adaptive sampling makes the required number of sampled columns/rows invariant to the distribution of the data.}
The proposed approach is amenable to online implementation and an online scheme is proposed.

\end{abstract}


\IEEEpeerreviewmaketitle

\section{Introduction}
\IEEEPARstart{S}{uppose}
we are given a data matrix $\bD \in \mathbb{R}^{N_1 \times N_2} $, which can be expressed as
\begin{eqnarray}
\bD= \bL+\bS,
\label{eq1}
\end{eqnarray}
where $\bL$ is a low rank (LR) matrix and $\bS$ a sparse matrix with arbitrary unknown support, whose entries can have arbitrarily
large magnitude. 
Many important applications in which the data under study can be naturally modeled using (\ref{eq1}) were discussed in \cite{lamport2}. The cutting-edge Principal Component Pursuit approach developed in \cite{lamport2 , lamport1},
directly decomposes $\bD$ into its LR and sparse components by solving the convex program 
%
\begin{eqnarray}
\begin{aligned}
 \underset{\dot{\bL},\dot{\bS}}{\min} \: \lambda\|\dot{\bS}\|_1  + \|\dot{\bL} \|_* \quad
  \text{subject to} \quad \dot{\bL} + \dot{\bS} = \bD \: \\
\end{aligned}
\label{eq2}
\end{eqnarray}
where $\|.\|_1$ is the $\ell_1$-norm,  $\|.\|_*$ is the nuclear norm  and $\lambda$  determines the trade-off between the sparse and LR components \cite{ lamport1}.
The convex program (\ref{eq2}) can precisely recover both the LR and sparse components if the columns and rows subspace of $\bL$ are sufficiently incoherent with the standard basis and the non-zero elements of $\bS$ are sufficiently diffused \cite{lamport1}.
Although the problem in (\ref{eq2}) is convex, its computational complexity is intolerable with large volumes of high-dimensional data. 
Even the efficient iterative algorithms proposed in \cite{lamport53, lamport23} have prohibitive computational and memory requirements 
in high-dimensional settings. 
%

\smallbreak
\noindent\textbf{Contributions:} This paper proposes a new randomized decomposition approach, which extracts the LR component in two consecutive steps. First, the column-space  (CS) of $\bL$ is learned from a small subset of the columns of the data matrix.
Second, the representation of the columns of $\bL$ with respect to the learned CS is obtained from a small subset of the rows.
 Unlike conventional decomposition that uses the entire data, we only utilize a small data sketch, and solve two low-dimensional optimization problems in lieu of one high-dimensional matrix decomposition problem (\ref{eq2}) resulting in significant running time speed-ups.

To the best of our knowledge, it is shown here for the first time that the sufficient number of randomly sampled columns/rows scales linearly with the rank $r$ and the coherency parameter of $\bL$ even with uniform random sampling.  Also, in contrast to the existing randomized approaches \cite{new1,lamport8}, which use blind uniform random sampling, we propose a new methodology for efficient column/row sampling. When the columns/rows of $\bL$ are not distributed uniformly in the CS/row space (RS) of $\bL$, which prevails much of the real world data, the proposed sampling approach is shown to achieve significant savings in data usage compared to uniform random sampling-based methods that require remarkable portions of the data.
 \textcolor{black}{The analysis presented shows that the proposed adaptive sampling procedure can make the required number of sampled columns/rows invariant to the data distribution. }
In addition, the proposed sampling algorithms can be independently used for feature selection from high-dimensional data.

In the presented approach, once the CS is learned, each column is decomposed efficiently and independently using the proposed randomized vector decomposition method. Unlike most existing approaches, which are batch-based, this unique feature enables applicability to online settings. 
The presented vector decomposition method can be independently used in many applications as an efficient vector decomposition algorithm or for efficient linear decoding \cite{lamport17}. 

\subsection{Notation and definitions}
We use bold-face upper-case letters to denote matrices and bold-face lower-case letters to denote vectors. Given a matrix $\bL$, $\|\bL\|$ denotes its spectral norm, $\|\bL\|_F $ its Frobenius norm,  and $\| \bL \|_\infty$ the infinity norm, which is equal to the maximum absolute value of its elements. In an $N$-dimensional space, $\be_i$ is the $i^{\text{th}}$ vector of the standard basis (i.e., the $i^{\text{th}}$ element of $\be_i$ is equal to one and all the other elements are equal to zero). The notation $A = [\:]$ denotes an empty matrix and the matrix $\bA = [ \bA_1 \: \bA_2 \: ... \: \bA_n] $
is the column-wise concatenation of the matrices $\{ \bA_i\}_{i=1}^n$. Random sampling refers to sampling without replacement.

\section{Background and Related Work}
\subsection{Exact LR plus sparse matrix decomposition}

The incoherence of the CS and RS of $\bL$ is an important requirement for the identifiability of the decompostion problem in (\ref{eq1}) \cite{lamport2 , lamport1}. For the LR matrix $\bL$ with rank $r$ and compact SVD $\bL = \bU \mathbf{\Sigma} \bV^T$ (where $\bU \in \mathbb{R}^{N_1 \times r}$, $\mathbf{\Sigma} \in \mathbb{R}^{r \times r}$ and $\bV \in \mathbb{R}^{N_2 \times r}$), the incoherence condition is typically defined through the requirements \cite{lamport1,lamport2}
\begin{eqnarray}
\begin{aligned}
& \underset{i}{\max} \|\bU^T\be_i\|_2^2 \leq \frac{\mu r}{N_1} \: , \: \underset{i}{\max} \|\bV^T\be_i\|_2^2 \leq \frac{\mu r}{N_2} \: \\
& \text{and} \: \: \|\bU\bV^T\|_{\infty} \leq \sqrt{\frac{\mu r}{N_2 N_1}}
\label{eq9}
\end{aligned}
\end{eqnarray}
for some parameter $\mu$ that bounds the projection of the standard basis $\{\be_i\}$ onto the CS and RS.  Other useful measures for the coherency of subspaces are given in \cite{lamport3} as,
\begin{eqnarray}
\gamma(\bU) \hspace{-0.5mm}= \hspace{-.5mm}\sqrt{N_1}
\max_{i,j} |\bU(i,j)|, \gamma(\bV) \hspace{-0.5mm}= \hspace{-0.5mm}\sqrt{N_2}
\max_{i,j} |\bV(i,j)|,
 \label{eq10}
\end{eqnarray}
where $\gamma(\bU)$ and $\gamma(\bV)$ bound the coherency of the CS and the RS, respectively. When some of the elements of the orthonormal basis of a subspace are too large, the subspace is coherent with the standard vectors. 

The decomposition of a data matrix into its LR and sparse components was analyzed in \cite{lamport1 , lamport2}, and sufficient conditions for exact recovery using the convex minimization (\ref{eq2}) were derived. In \cite{lamport2}, the sparsity pattern of the sparse matrix is selected uniformly at random following the so-called Bernoulli model to ensure that the sparse matrix is not LR with overwhelming probability. In this model, which is also used in this paper, each element of the sparse matrix can be non-zero independently with a constant probability. 
Without loss of generality (w.l.o.g.), suppose that $N_2 \leq N_1$. The following lemma states the main result of \cite{lamport2}.
\begin{lemma}[\text{Adapted from \cite{lamport2}}]
Suppose that the support set of $\bS$ follows the Bernoulli model with parameter $\rho$. The convex program (\ref{eq2}) with $\lambda=\frac{1}{\sqrt{N_1}}$ yields the exact decomposition with probability at least $1-c_1 {N_1}^{-10}$  provided that
\begin{eqnarray}
r \leq \rho_r N_2 {\mu}^{-1} \left(\log(N_1)\right)^{-2} \quad , \quad \rho \leq {\rho}_s
\label{eq12}
\end{eqnarray}
where $\rho_s$, $c_1$ and $\rho_r$ are numerical constants.
\label{lm1}
\end{lemma}
The optimization problem in (\ref{eq2}) is convex and can be solved using standard techniques such as interior point methods \cite{lamport1}. Although these methods have fast convergence rates, their usage is limited to small-size problems due to the high complexity of computing a step direction. Similar to the iterative shrinking algorithms for $\ell_1$-norm and nuclear norm minimization, a family of iterative algorithms for solving the optimization problem (\ref{eq2}) were proposed in \cite{lamport23,lamport53}. However, they also 
require working with the entire data. For example, the algorithm in \cite{lamport23} requires computing the Singular Value Decomposition (SVD) of an $N_1 \times N_2 $ matrix in every iteration.

\subsection{Randomized approaches}
%
Owing to their inherent low-dimensional structures, robust principal component analysis (PCA) and low rank plus sparse matrix decomposition can be conceivably solved using small data sketches, i.e., a small set of random observations of the data \cite{lamport16,new3,feldman2010coresets,lamport8,lamport14,rahmani2015randomized,rahmani2016subspace}. 
In \cite{lamport16}, it was shown based on simple degree-of-freedom analysis that the LR and sparse components can be precisely recovered using a small set of random linear measurements of $\bD$. A convex program was proposed in \cite{lamport16} to recover these components using random matrix embedding with a polylogarithmic penalty factor in sample complexity, albeit the formulation also requires solving a high-dimensional optimization problem.

The iterative algorithms which solve (\ref{eq2}) have complexity $\calO(N_1 N_2 r)$ per iteration since they compute the partial SVD of $N_1 \times N_2$ matrices \cite{lamport23}. To reduce complexity, GoDec \cite{zhou2011godec} uses a randomized method to efficiently compute the SVD, and the decomposition algorithm in \cite{mu2011accelerated} minimizes the rank of $\mathbf{\Phi} \bL$ instead of $\bL$, where $\mathbf{\Phi}$ is a random projection matrix. However, these approaches do not have provable performance guarantees and their memory requirements scale with the full data dimensions.
Another limitation of the algorithm in \cite{mu2011accelerated} is its instability since different
random projections may yield different results. 

The divide-and-conquer approach in \cite{new1} (and a similar algorithm in \cite{liu2011solving}),  
can achieve super-linear speedups over full-scale
matrix decomposition. This approach forms an estimate of $\bL$ by combining two low-rank approximations obtained from submatrices formed from sampled rows and columns of $\bD$ using the generalized Nystr\"om method \cite{goreinov1997theory}.
%
Our approach also achieves super-linear speedups in decomposition, yet is fundamentally different from \cite{new1} and offers several advantages. 
First, our approach is a \emph{subspace-pursuit approach} that focuses on subspace learning in a structure-preserving data sketch. Once the CS is learned, each column of the data is decomposed independently using a proposed randomized vector decomposition algorithm. Second, unlike \cite{new1}, which is a batch approach that requires to store the entire data, the structure of the proposed approach naturally lends itself to online implementation (c.f. Section \ref{sec:online}), which could be very beneficial when the data comes in on the fly \cite{mardani2015subspace,mardani2015online,zhan2014performance,guo2014online}. Third,
while the analysis provided in \cite{new1} requires roughly $\calO (r^2 \mu^2 \max (N_1 , N_2))$ random observations to ensure exact decomposition with high probability (whp), we show that the order of sufficient number of random observations depends linearly on the rank and the coherency parameter even if uniform random sampling is used.  \textcolor{black}{ Fourth, the structure of the proposed approach enables us to leverage adaptive sampling strategies for challenging and realistic scenarios in which the columns and rows of $\bL$ are not uniformly distributed in their respective subspaces, or when the data exhibits additional structures (e.g. clustering) (c.f. Sections \ref{sec:eff_sampling},\ref{sec:alg_eff_sampling}). It is shown that the proposed adaptive sampling scheme can make the number of randomly sampled columns/rows required for exact decomposition invariant to the data distribution. In such settings, the uniform random sampling used in \cite{new1} requires significantly larger amounts of data to carry out the decomposition.}

\section{Structure of the Proposed Approach and Theoretical Result}
%
%




%
In this section, the structure of the proposed randomized decomposition method is presented. A step-by-step analysis of the proposed approach is provided and sufficient conditions for exact decomposition are derived. Theorem \ref{thm:main_result} stating the main theoretical result of the paper is presented at the end of this section. The proofs of the lemmas and the theorem are deferred to the appendix.

Let us rewrite (\ref{eq1}) as
$
\bD=\bU\bQ+\bS,
$
where $\bQ=\boldsymbol{\Sigma} \bV $. 
The representation matrix $\bQ \in \mathbb{R}^{r \times N_2} $ is a full row rank matrix that contains the expansion of the columns of $\bL$ in the orthonormal basis $\bU$. 
The first step of the proposed 
approach aims to learn the CS of $\bL$ using a subset of the columns of $\bD$, and in the second step the representation matrix is obtained using a subset of the rows of $\bD$.

Let $\:\calU$ denote the CS of $\bL$. Fundamentally, $\calU$ can be obtained from a small subset of the columns of $\bL$. However, since we do not have direct access to the LR matrix, a random subset of the columns of  $\bD$ is first selected. Hence, the matrix of sampled columns $\bD_{s1}$ can be written as
$\bD_{s1}=\bD \bS_{1}$,
where $\bS_1 \in \mathbb{R}^{N_2 \times m_1} $ is the column sampling matrix and $m_1$ is the number of  selected columns. The matrix of selected columns can be written as
\begin{eqnarray}
\bD_{s1}=\bL_{s1} + \bS_{s1},
\label{eq19}
\end{eqnarray}
\textcolor{black}{where $\bL_{s1}$ and $\bS_{s1}$ are its LR and sparse components, respectively.
The idea is to decompose the sketch $\bD_{s1}$ into its LR and sparse components to learn the CS of $\bL$ from the CS of $\bL_{s1}$. 
Note that the columns of $\bL_{s1}$ are a subset of the columns of $\bL$ since $\bL_{s1} = \bL \bS_1$. Should we be able to decompose $\bD_{s1}$ into its exact LR and sparse components (c.f. Lemma \ref{lm4}), we also need to ensure that the columns of $\bL_{s1}$ span $\calU$. The following lemma establishes that a small subset of the columns of $\bD$ sampled uniformly at random contains sufficient information (i.e., the columns of the LR component of the sampled data span $\calU$) if the RS is incoherent. 
}
\begin{lemma}
Suppose $m_1$ columns are sampled uniformly at random from the matrix $\bL$ with rank $r$. If
\begin{eqnarray}
m_1 \ge r\gamma^2  (\bV) \max \left( c_2 \log r , c_3 \log\left(\frac{3}{\delta}\right) \right),
\label{eq27}
\end{eqnarray}
the selected columns of the matrix $\bL$ span the CS of $\bL$ with probability at least $1-\delta$ for constants $c_2$ and $c_3$.
\label{lm3}
\end{lemma}
Thus, if $\gamma (\bV)$ is small (i.e., the RS is not coherent), a small set of randomly sampled  columns can span $\calU$. Based on Lemma \ref{lm3}, if $m_1$ satisfies (\ref{eq27}), then $\bL$ and $\bL_{s1}$ will have the same CS whp. The following optimization problem (of dimensionality $N_1 m_1$) is solved to decompose $\bD_{s1}$ into its LR and sparse components.
\begin{eqnarray}
\begin{aligned}
& \underset{\dot{\bL}_{s1},\dot{\bS}_{s1}}{\min} \: \frac{1}{\sqrt{N_1}}\| \dot{\bS}_{s1} \|_1  + \| \dot{\bL}_{s1} \|_* \quad\\
& \text{subject to} \quad  \dot{\bL}_{s1} + \dot{\bS}_{s1} = \bD_{s1} . \\
\end{aligned}
\label{eq20}
\end{eqnarray}
Thus, the columns subspace of the LR matrix can be recovered by finding the columns subspace of $ \bL_{s1} $. Our next lemma establishes that (\ref{eq20}) yields the exact decomposition using roughly $m_1 = \calO (\mu r)$ randomly sampled columns. To simplify the analysis, in the following lemma it is assumed that the CS of the LR
matrix is sampled from the random orthogonal model \cite{lamport4}, i.e., the columns of $\bU$ are selected
uniformly at random among all families of $r$-orthonormal vectors.
\begin{lemma}
Suppose the columns subspace of $\bL$ is sampled from the random orthogonal model, $\bL_{s1}$ has the same column subspace of $\bL$ and the support set of $\bS$ follows the Bernoulli model with parameter $\rho$. In addition, assume that the columns of $\bD_{s1}$ were sampled uniformly at random. If
\begin{eqnarray}
m_1  \ge {\frac{r}{\rho}_r} \mu^{'} (\log N_1)^2 \: \: \: \text{and} \: \: \: \rho \leq \rho_s \:,
\label{eq28}
\end{eqnarray}
then (\ref{eq20}) yields the exact decomposition with probability at least $1-c_8 N_1^{-3}$, where
\begin{eqnarray}
\mu^{'} \hspace{-1.5mm} = \hspace{-.5mm} \max \hspace{-.75mm}\left( \hspace{-.75mm}\frac{c_7 \max (r,\log N_1)}{r} , 6\gamma^2(\bV) , (c_9\gamma(\bV) \log N_1)^2 \hspace{-.75mm}\right)
\label{eq30}
\end{eqnarray}
and $c_7$, $c_8$ and $c_9$ are constant numbers provided that $N_1$ is greater than the RHS of the first inequality of (\ref{eq28}).
\label{lm4}
\end{lemma}
Therefore, according to Lemma \ref{lm4} and Lemma \ref{lm3}, the CS of $\bL$ can be obtained using roughly $O(r \mu)$ uniformly sampled data columns. Note that $m_1 \ll N_1$ for high-dimensional data as $m_1$ scales linearly with $r$. Hence, the requirement that $N_1$ is also greater than the RHS of the first inequality of (\ref{eq28}) is by no means restrictive and is naturally satisfied.

Define $\hat{\bU}$ as an orthonormal basis for the learned CS.
 An arbitrary column $\bd_i$ of $\bD$ can be written as $\bd_i = \bU \bq_i + \bs_i$, where $\bq_i$ and $\bs_i$ are the corresponding columns of $\bQ$ and $\bS$, respectively. Thus, $\bd_i - \bU \bq_i$ is a sparse vector. This suggests that $\bq_i$ can be learned using the minimization
\begin{eqnarray}
\underset{\hat{\bq_i}}{\min} \|\bd_i-\hat{\bU} \hat{\bq_i}\|_1  \:,
\label{firstrow}
\end{eqnarray}
where the $\ell_1$-norm is used as a surrogate for the $\ell_0$-norm to promote a sparse solution \cite{lamport17 , lamport3}.
The optimization problem (\ref{firstrow}) is similar to a system of linear equations with $r$ unknown variables and $N_1$ equations. Since $r \ll N_1$, the idea is to learn $\bq_i$ using only a small subset of the equations. Thus, we propose the following vector decomposition program
\begin{eqnarray}
\underset{\hat{\bq_i}}{\min} \|\bS_2^T \bd_i- \bS_2^T \hat{\bU} \hat{\bq_i}\|_1  \:,
\label{n_eq1}
\label{secondrows}
\end{eqnarray}
where $\bS_2 \in \mathbb{R}^{N_1 \times m_2} $ selects $m_2$ rows of $\hat{\bU}$ (and the corresponding $m_2$ elements of $\bd_i$).

First, we have to ensure that the rank of $\bS_2^T \bU$ is equal to the rank of $\bU$, for if $\bq^{*}$ is the optimal point of (\ref{secondrows}), then $\bU \bq^{*}$ will be the LR component of $\bd_i$. According to Lemma \ref{lm3}, $m_2 = \calO (r \gamma^2(\bU))$, is sufficient to preserve the rank of $\bU$ when the rows are sampled uniformly at random. 
%
In addition, the following lemma establishes that if the rank of $\bU$ is equal to the rank of $\bS_2^T \bU$, then the sufficient value of $m_2$ for (\ref{n_eq1}) to yield the correct columns of $\bQ$ whp is linear in $r$.

\begin{lemma}
Suppose that the rank of $\bS_2^T \bU$ is equal to the rank of $\bL$ and assume that the CS of $\bL$ is sampled from the random orthogonal model. The optimal point of (\ref{secondrows}) is equal to $\bq_i$ with probability at least $(1-3\delta^{'})$ provided that
\begin{eqnarray}
\begin{aligned}
&\rho \leq \frac{0.5}{r \beta \left( c_6  \kappa \log \frac{N_1 }{\delta^{'}} + 1 \right)}~, \\
&m_2 \ge \max \bigg( \frac{ 2 r \beta (\beta - 2) \log \left( \frac{ 1 }{\delta^{'}} \right)}{3(\beta - 1)^2} \left( c_6  \kappa \log \frac{N_1 }{\delta^{'}} +1 \right) ,  \\
&\qquad\qquad\qquad\qquad\qquad\qquad  c_5 (\log \frac{N_1  }{\delta^{'}})^2 , \sqrt[6]{\frac{3}{\delta^{'}}} \bigg)
\end{aligned}
\label{eq38}
\end{eqnarray}
where $\kappa= \frac{\log N_1}{r}$, $c_5$ and $c_6$ are constant numbers and $\beta$ can be any real number greater than one.
\label{lm8}
\end{lemma}
Therefore, we can obtain the LR component of each column using a random subset of its elements. Since (\ref{firstrow}) is an $\ell_1$-norm minimization, we can write the representation matrix learning problem as
\begin{eqnarray}
\underset{\dot{\bQ}}{\min} \| \bD_{s2} - \hat{\bU}_{s2} \dot{\bQ}\|_1,
\label{eq23}
\end{eqnarray}
where $\hat{\bU}_{s2} = \bS_2^T \hat{\bU}$.
Thus, (\ref{eq23}) learns $\bQ$ using a  subset of the rows of $\bD$ as $\bS_2^T \bD$ is the matrix formed from $m_2$ sampled rows of $\bD$.

As such, we solve two low-dimensional subspace pursuit problems (\ref{eq20}) and (\ref{eq23}) of dimensions $N_1 m_1$ and $N_2 m_2$, respectively, instead of an $N_1 N_2$-dimensional decomposition problem (\ref{eq2}), and use a small random subset of the data to learn $\bU$ and $\bQ$.
The table of Algorithm 1 explains the structure of the proposed approach.



We can readily state the following theorem which establishes sufficient conditions for Algorithm 1 to yield exact decomposition. 

\begin{theorem}
\label{thm:main_result}
\textcolor{black}{Suppose the CS of the LR matrix is sampled from the random orthogonal model and the support set of $\bS$ follows the Bernoulli model with parameter $\rho$. Also, it is assumed that Algorithm 1 samples the columns and rows uniformly at random.  If for any small $\delta > 0$, $m_1$ satisfies the inequalities (\ref{eq27}) and (\ref{eq28}), $\rho$ satisfies inequality (\ref{eq38}) with $\delta^{'} = \delta/N_2$, $m_2$ satisfies inequality (\ref{eq38}) with $\delta^{'} = \delta/N_2$ and
also
\begin{eqnarray}
\begin{aligned}
m_2 & \ge  r \log N_1 \max \big( c_2^{'} \log r , c_3^{'} \log \frac{3}{\delta} \big) \quad , \quad \rho \leq \rho_s
\end{aligned}
\label{eq25}
\end{eqnarray}
where $\{c_i \}_{i = 1}^9$, $c_2^{'}$ and $c_3^{'}$ are constant numbers, $\mu^{'}$ is equal to (\ref{eq30}), $\kappa= \frac{\log N_1}{r}$, and $\beta$ can be any real number greater than one,
then the proposed approach (Algorithm 1) yields exact decomposition with probability at least $(1-5\delta - 3r N_1^{-7}  - c_8 N_1^{-3})$ provided that $N_1$ is greater than the RHS of the first inequality of (\ref{eq28}).
}
\end{theorem}

Theorem \ref{thm:main_result} guarantees that the LR component can be obtained using a small subset of the data.
The randomized approach has two main advantages. First, it significantly reduces the memory/storage requirements since it only uses a small data sketch and solves two low-dimensional optimization problems versus one large problem. Second, the proposed approach has $\calO(\max (N_1 , N_2)  \times \max(m_1 , m_2) \times r)$ per-iteration running time complexity, which is significantly lower than $\calO(N_1 N_2 r)$ per iteration for full scale decomposition (\ref{eq2}) \cite{lamport53 , lamport23} implying remarkable speedups for big data. For instance, consider $\bU$ and $\bQ$ sampled from $\calN (0 , 1)$, $r = 5$, and $\bS$ following the Bernoulli model with $\rho = 0.02$. For values of $N_1 = N_2$ equal to 500, 1000, 5000, $10^4$ and  $2\times 10^4$, if $m_1 = m_2 = 10r$, the proposed approach yields the correct decomposition with 90, 300, 680, 1520 and 4800 - fold speedup, respectively, over directly solving (\ref{eq2}).
 %
%
\begin{algorithm}
\caption{Structure of Proposed Approach}
{\footnotesize
\textbf{Input}: Data matrix $\bD \in \mathbb{R}^{N_1 \times N_2} $\\
\textbf{1. Initialization}: Form  column sampling matrix $\bS_1 \in \mathbb{R}^{N_2 \times m_1} $ and row sampling matrix $\bS_2 \in \mathbb{R}^{N_1 \times m_2} $.\\
\textbf{2. CS Learning}\\
\textbf{2.1} Column sampling: Matrix $\bS_1$ samples $m_1$ columns of the given data matrix, $\bD_{s1} = \bD \bS_1$. \\
\textcolor{black}{
\textbf{2.2} CS learning: Matrix $\hat{\bL}_{s1}$ is obtained as the LR component of $\bD_{s1}$  (\ref{eq20}). \\
\textbf{2.3} CS calculation: Matrix $\hat{\bU}$ is formed as an orthonormal basis for the CS of $\hat{\bL}_{s1}$. \\
\textbf{3. Representation Matrix Learning}\\
\textbf{3.1} Row sampling: Matrix $\bS_2$ samples $m_2$ rows of the given data matrix, $\bD_{s2} = \bS_2^T \bD$.\\
\textbf{3.2} Matrix $\hat{\bQ}$ is obtained as the optimal point of
(\ref{eq23}).}
\\
\textbf{Output:} If $\hat{\bU}$ is an orthonormal basis for the learned CS and $\hat{\bQ}$ is the obtained representation matrix, then $\hat{\bL} = \hat{\bU} \hat{\bQ}$ is the obtained LR component.
}
\end{algorithm}
%
%

\section{Efficient Column/Row Sampling}
In sharp contrast to randomized algorithms for matrix approximations rooted in numerical linear algebra (NLA) \cite{ghadim1,mahoney2011randomized}, which seek to compute matrix approximations from sampled data using importance sampling, in matrix decomposition and robust PCA we do not have direct access to the LR matrix to measure how informative particular columns/rows are. As such, the existing randomized algorithms for matrix decomposition and robust PCA\cite{new1 , new3 , lamport8} have predominantly relied upon uniform random sampling of columns/rows.  

In Section \ref{sec:non_uniform}, we briefly describe the implications of non-uniform data distribution and show that uniform random sampling may not be favorable for data matrices exhibiting some structures that prevail much of the real datasets. In Section \ref{sec:eff_sampling}, we demonstrate an efficient column sampling strategy which will be integrated with the proposed decomposition method. The decomposition method with efficient column/row sampling is presented in Section \ref{sec:alg_eff_sampling}. 

\subsection{Non-uniform data distribution}
\label{sec:non_uniform}
\begin{figure}
	\centering
    \includegraphics[width=0.41\textwidth]{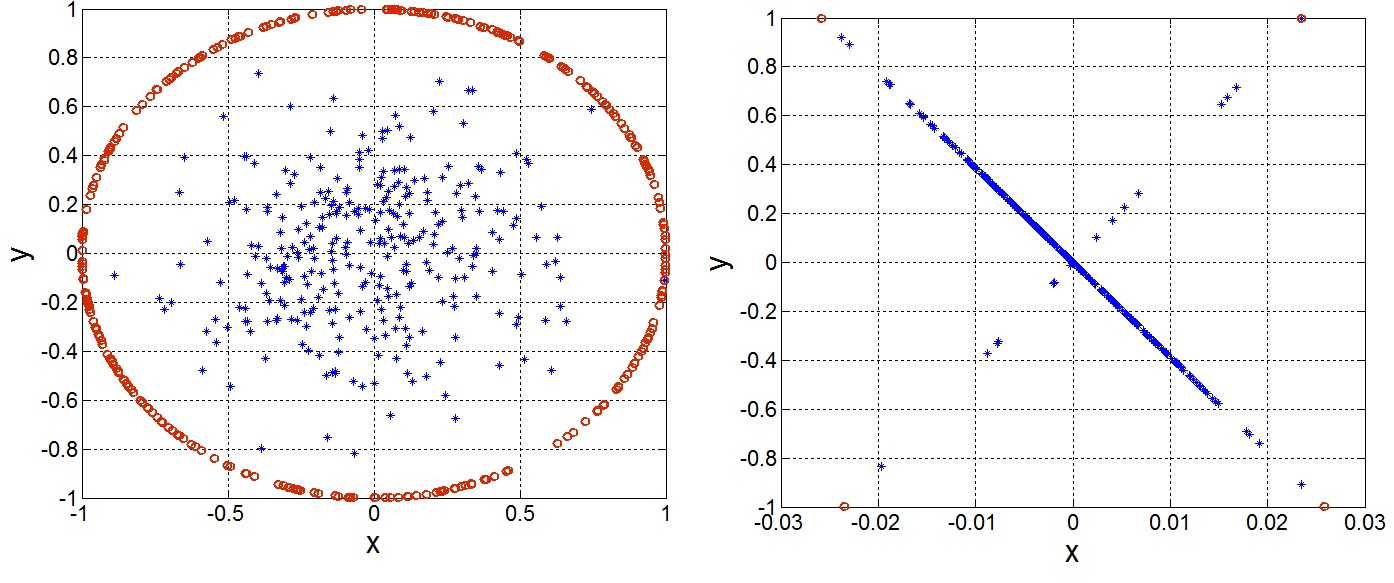}
    \vspace{-0.3cm}
    \caption{Data distributions in a two-dimensional subspace. The red points are the normalized data points.}
    \label{fig:innovationfig}
\end{figure}

When data points lie in a low-dimensional subspace, a small subset of the points can span the subspace. However, uniform random sampling is only effective when the data points are distributed uniformly in the subspace. To clarify, Fig.~\ref{fig:innovationfig} shows two scenarios for a set of data points in a two-dimensional subspace. In the left plot, the data points are distributed uniformly at random. 
In this case, two randomly sampled data points can span the subspace whp. In the right plot, 95 percent of the data lie on a one-dimensional subspace, thus we may not be able to capture the two-dimensional subspace from a small random subset of the data points.
\begin{figure}
	\centering
    \includegraphics[width=0.38\textwidth]{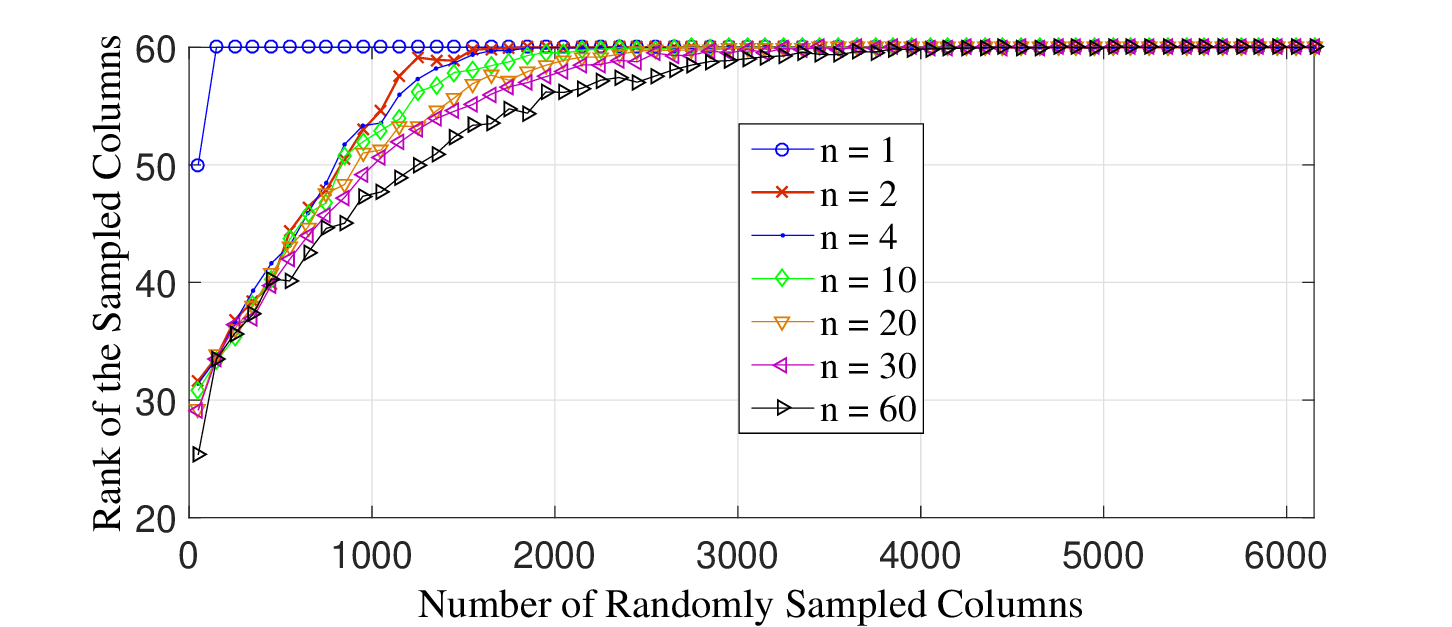}
    \vspace{-0.25cm}
    \caption{The rank of a set of uniformly sampled columns for different number of clusters. }
        \label{fig:mumv}
\end{figure}

In practice, the data points in a low-dimensional subspace may not be uniformly distributed, but rather exhibit some additional structures.
A prevailing structure in many modern applications is clustered data \cite{new5}. For example, user ratings for certain products (e.g. movies) in recommender systems are not only LR due to their inherent correlations, but also exhibit additional clustering structures owing to the similarity of the preferences of individuals from similar backgrounds (e.g. education, culture, or gender) \cite{new5}.



To further show that uniform random sampling falls short when the data points are not distributed uniformly in the subspace, consider a
matrix $\bG \in \mathbb{R}^{2000 \times 6150}$ generated as
$
\bG = [ \bG_1 \: \bG_2 \: ... \: \bG_n ]
$. For $1 \leq i \leq \frac{n}{2}$, $\bG_i = \bU_i \bQ_i\: ,$
where $\bU_i \in \mathbb{R}^{2000 \times \frac{r}{n}}$, $\bQ_i \in \mathbb{R}^{\frac{r}{n} \times \frac{200 r}{n}}$. For $n/2+1 \leq i \leq n$, $\bG_i = \bU_i \bQ_i\: ,$
where $\bU_i \in \mathbb{R}^{2000 \times \frac{r}{n}}$, $\bQ_i \in \mathbb{R}^{\frac{r}{n} \times \frac{5 r}{n}}$.
The elements of $\bU_i$ and $\bQ_i$ are sampled independently from a normal $\calN(0 , 1)$ distribution.
The parameter $r$ is set equal to 60, thus the rank of $\bG$ is equal to 60 whp. Fig.~\ref{fig:mumv} illustrates the rank of the randomly sampled columns versus the number of sampled columns for different number of clusters $n$.  When $n=60$, it turns out that we need to sample more than half of the columns to span the CS, 
so we cannot evade high-dimensionality with uniform random column/row sampling.
\textcolor{black}{
Thus, when the distribution of the clustered data is less uniform, more randomly sampled columns are needed  to capture the CS.}
This is indeed confirmed by Lemma \ref{lm3}, which established that the sufficient number of randomly sampled columns is proportional to the RS coherency.
 In the following lemmas, it is shown that the RS coherency increases if the distribution of the columns is less uniform. Lemma \ref{lm:coherece_study} provides an upper bound on the coherency of the RS and Lemma \ref{Lm:lowerbound} confirms that the upper bound is tight by establishing a converse.
In these lemmas, it is assumed that the columns lie in a union of linear subspaces as per the following assumption.
\begin{assumption}
The matrix $\bL$ can be represented as $\bL = [\bU_1 \bQ_1 \: ... \: \bU_n \bQ_n]$. The CS of $\{ \bU_i \in \mathbb{R}^{N_1 \times r/n} \}_{i =1}^{n}$ are  random $r/n$-dimensional subspaces in $\mathbb{R}^{N_1}$. The RS of
$\{ \bQ_i \in \mathbb{R}^{r/n \times n_i} \}_{i =1}^{n}$ are random $r/n$-dimensional subspaces in $\{ \mathbb{R}^{n_i} \}_{i=1}^n$, respectively, $\sum_{i = 1}^n n_i = N_2$, and $\underset{i}{\min}\: {n_i} \gg r/n$.
\end{assumption}


\begin{lemma}
Suppose $\bL$ follows Assumption 1. If the rank of $\bL$ is equal to $r$, then
\begin{eqnarray}
\begin{aligned}
 & \mathbb{P} \Bigg[ \underset{i}{\max} \: \|\be_i^T \bV \|_2^2  > \frac{c_7 r \varphi_1 }{ N_2 } \left(  \frac{1}{n}  \frac{N_2}{ \underset{k}{\min} \: n_k } \right)  \Bigg] \leq 2  \sum_{i=1}^n n_i^{-3} \\
 & \mathbb{P} \left[
\gamma^2(\bV)  >  20 \log (\underset{i}{\max}\: n_i) \frac{ N_2}{ \underset{i}{\min}\: n_i } \right] \leq \frac{3r}{n} \sum_{i=1}^n n_i^{-7}
\end{aligned}
\end{eqnarray}
where $\varphi_1 = \frac{\max (r/n , \log \underset{k}{\max} \: n_k)}{r/n} \:.$
\label{lm:coherece_study}
\end{lemma}
\textcolor{black}{
\begin{lemma}
If $\bL$ follows Assumption 1, the rank of $\bL$ is equal to $r$, $r/n \ge 18 \log \underset{i}{\max} \: n_i$ and $ n_i \ge 96 \frac{r}{n} \log n_i, \: 1 \leq i \leq n$, then
\begin{eqnarray}
\mathbb{P} \left[ \underset{i}{\max} \: \| \bV^T \be_i \|_2^2  <  \frac{0.5 \: r}{N_2} \left( \frac{1}{n} \frac{N_2}{\underset{i}{min} \: n_i} \right) \right] \leq 2 \sum_{i = 1}^n n_i^{- 5}  \:.
\end{eqnarray}
\label{Lm:lowerbound}
\end{lemma}
}
Based on Lemma \ref{lm:coherece_study}, the RS coherency of $\bL$ is linear with $\frac{N_2}{\underset{i}{\min}\: n_i}$, confirming that the RS coherency increases if the distribution of the columns within the CS is less uniform.


\subsection{Efficient column sampling method}
\label{sec:eff_sampling}
Column sampling is 
widely used for dimensionality reduction and feature selection \cite{lamport14,new4}.
In the column sampling problem, the LR matrix (or the matrix whose span is to be approximated with a small set of its columns) is available. Thus, the columns are sampled based on their importance, measured by the so-called leverage scores \cite{ghadim1}, as opposed to blind uniform sampling.   
We refer the reader to \cite{lamport14 , new4} and references therein for more information about efficient column sampling methods. 

\textcolor{black}{Next, we present a sampling approach to be used in Section \ref{sec:alg_eff_sampling} where the proposed decomposition algorithm with efficient sampling is presented. The proposed sampling strategy is inspired by the approach in \cite{new4} in the context of volume sampling.
%
Algorithm 2 details the presented sampling procedure. Given a matrix $\bA$ with rank $r_A$, the algorithm aims to sample a small subset of columns that span its CS. The first column is sampled uniformly at random or based on a judiciously chosen probability distribution \cite{ghadim1}. The next columns are selected sequentially so as to maximize the novelty to the span of the selected columns. As shown in step 2.2 of Algorithm 2, a design threshold $\tau$ is used to decide whether a given column brings sufficient novelty to the sampled columns by thresholding the $\ell_2$-norm of its projection on the complement of the span of the sampled columns. The threshold $\tau$ is set to zero in a noise-free setting. Once the selected columns are believed to span the CS of $\bA$, they are removed from $\bA$. This procedure is repeated $C$ times (using the remaining columns). In each time, the algorithm finds $r_A$ columns spanning the CS of $\bA$. After every iteration, the rank of the matrix of remaining columns is bounded above by $r_A$. As such, the algorithm samples approximately $m_1 \approx C r_A$ columns in total.
 In the proposed decomposition method with efficient column/row sampling (c.f. Sec. \ref{sec:alg_eff_sampling}), we set $C$ large enough to ensure that the selected columns form a low rank matrix. }

\vspace{-.1cm}
\begin{algorithm}
\caption{Efficient Sampling from LR Matrices}
{\footnotesize
\textbf{Input:} Matrix $\bA$. \\
\textbf{1. Initialize}\\
\textbf{1.1} The parameter $C$ is chosen as an integer greater than or equal to one. The algorithm finds $C$ sets of linearly dependent columns. \\
\textbf{1.2} Set $\calI = \emptyset $ as the index set of the sampled columns and set $v = \tau$, $\bB = \bA$ and $\bC = [\:]$.   \\
 \textcolor{black}{ \textbf{2. Repeat} C Times}
\\
\textbf{2.1} Let $\bb$ be a non-zero randomly sampled column from $\bB$ with index $i_b$. Update $\bC$ and $\calI$ as $\bC = [\bC \: \: \bb]$, $\calI = \{ \calI \: , \: i_b \}$.  \\
\textbf{2.2} \textbf{While $v \ge \tau$}\\
\textbf{2.2.1} Set $\bE = \bP_c \bB \: , $
where $\bP_c$ is the projection matrix onto the complement space of $\spn(\bC)$. \\
\textbf{2.2.2} Define $\mathbf{f}$ as the column of $\bE$ with the maximum $\ell_2$-norm with index $i_f$. Update $\bC$, $\calI$ and $v$ as $\bC = [\bC \: \:  \: \mathbf{f}] \: \: , \: \: \calI = \{ \calI \: , \: i_f \} \: \: \text{and} \:\: v = \| \mathbf{f} \|_2 \: . $\\
\textbf{2.2 End While}\\
\textbf{2.3} Set $\bC = [\:]$ and set $\bB$ equal to $\bA$ with the columns indexed by $\calI$ set to zero. \\
 \textcolor{black}{ \textbf{2. End Repeat}}
\\
\textbf{Output:} The set $\calI$ contains the indices of the selected columns.
}
\end{algorithm}
%
\begin{figure}[t!]
	\centering
    \includegraphics[width=0.44\textwidth]{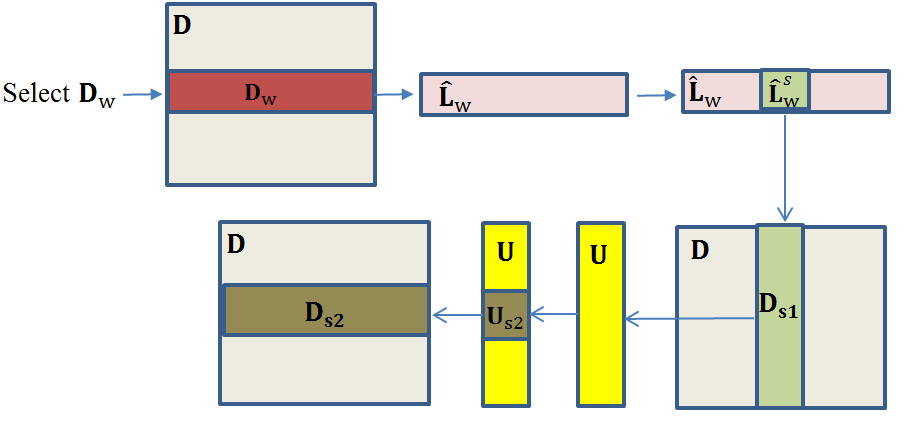}
    \vspace{-0.3cm}
    \caption{Visualization of the matrices defined in Section \ref{sec:alg_eff_sampling}. Matrix $\bD_w$ is selected randomly or using Algorithm 3 described in Section \ref{sec:both_col_row}.}
    \label{fig: visula 1}
\end{figure}

\begin{algorithm}
\caption{Decomposition with Adaptive Sampling }
{\footnotesize
\smallbreak
\textbf{1. Informative Column Sampling}\\
\textbf{1.1} Sample $C_r \hat{r}$ rows of $\bD$ uniformly at random to form $\bD_w$ where $\hat{r}$ is a known upper bound on $r$.\\
\textbf{1.2} Obtain $\hat{\bL}_w$  via (\ref{eq2}) as the LR component of $\bD_w$.
\\
\textbf{1.3} Apply Algorithm 2 to $\hat{\bL}_w$. \\
\textbf{1.4} Form the matrix $\bD_{s1}$ from the columns of $\bD$ corresponding to the sampled columns of $\hat{\bL}_w$.

\smallbreak
\textbf{2. CS Learning}\\
\textbf{2.1} The convex program (\ref{eq20}) is applied to  $\bD_{s1}$. \\
\textbf{2.2} Obtain $\hat{\bU}$ as an orthonormal basis for the CS of the calculated LR component of $\bD_{s1}$.

\smallbreak
\textbf{3. Representation Learning}\\
\textbf{3.1} Apply Algorithm 2 to $\hat{\bU}^T$  and define $\hat{\bU}_{s2}$ as the matrix of sampled rows.
\\
\textbf{3.2} Form the matrix $\bD_{s2}$ as the rows of $\bD$ corresponding to the sampled rows of $\hat{\bU}$.
\\
\textbf{3.3} Define $\hat{\bQ}$ as the optimal point of (\ref{eq23}).

\smallbreak
\textbf{Output:} The matrix $\hat{\bU}$ is the obtained basis for the column space of $\bL$ and $\hat{\bL} = \hat{\bU} \hat{\bQ}$ is the obtained LR matrix.

}
\end{algorithm}
\vspace{-.5cm}

\section{Proposed decomposition algorithm with efficient sampling}
\label{sec:alg_eff_sampling}
In this section, we present a modified decomposition algorithm (Algorithm 3) that replaces uniform random sampling with the efficient column/row sampling method (Algorithm 2). \textcolor{black}{We subsequently provide an analysis showing that the required number of sampled columns using the proposed approach can be invariant to the distribution of the columns -- even if their distribution is highly non-uniform}. The table of Algorithm 3 details the proposed method with efficient sampling.  In this algorithm, it is assumed that the \emph{rows} of $\bL$ are well distributed, in the sense that they do not align along any specific directions, such that $C_r r$ rows of $\bL$ sampled uniformly at random span its RS whp, for some constant $C_r$. In Section \ref{sec:both_col_row}, we dispense with this assumption. 
The proposed decomposition algorithm consists of three steps detailed next.  \textcolor{black}{ We support the idea of each step with some examples and theoretical analysis.}
\\
\subsubsection{Informative column sampling}
 The first important idea underlying the proposed sampling approach is to start sampling along the dimension that has the better distribution. For instance,
 consider an extreme scenario where only two columns of $\bG \in \mathbb{R}^{1000 \times 1000}$ are non-zero. In this case, with random sampling we need to sample almost all the columns to ensure that the sampled columns span the CS of $\bG$.
But, if the non-zero columns are non-sparse, a small subset
of randomly chosen rows of $\bG$ will span its row space. \textcolor{black}{As another example, suppose that the distribution of the columns of $\bL$ follows Assumption 1, i.e., the distribution of the columns admits a clustering structure. The following lemmas compare the sufficient numbers of randomly sampled columns and rows to capture the CS and RS, respectively. }
\textcolor{black}{
\begin{lemma}
Suppose $\bL$ follows Assumption 1. If $m_1$ columns of $\bL$ are sampled uniformly at random with replacement, the rank of $\bL$ is equal to $r$, and
\begin{eqnarray}
\begin{aligned}
m_1 \ge \left( 2 + \frac{3}{\xi_{\min}} \log\frac{2 n}{\delta} \right)  \frac{\xi_{\max} N_2}{ \underset{i}{\min}\: n_i }
\end{aligned}
\label{eq:suf_sample_c}
\end{eqnarray}
where
\begin{eqnarray}
\begin{aligned}
& \xi_{\min} = {10} \: c_7 \: {\max(r/n , \log \underset{i}{\min}\: n_i )} \log \frac{2 r}{\delta} \\
&\xi_{\max} = {10} \: c_7 \: {\max(r/n , \log \underset{i}{\max}\: n_i )} \log \frac{2 r}{\delta} \:,
\end{aligned}
\end{eqnarray}
then the sampled columns span the column space of $\bL$ with probability at least $1 - 2 \delta - 2 \sum_{i = 1}^n n_i^{-3}$.
\label{lm:columnspmle}
\end{lemma}
\begin{lemma}
Suppose $\bL$ follows Assumption 1 and $m_2$ rows of $\bL$ are sampled uniformly at random with replacement. If the rank of $\bL$
 is equal to $r$ and
\begin{eqnarray}
m_2 \ge 10 \: c_7 \: r \varphi_2 \log \frac{2 r}{\delta},
\label{eq:suff_row_com}
\end{eqnarray}
then the sampled rows span the row space of ~$\bL$ with probability at least $1 - \delta - 2 N_1^{-3}$, where $\varphi_2 = \frac{\max (r ,  \log N_1 )}{r}$.
\label{lm:row_compare}
\end{lemma}
The sufficient number of sampled rows indicated in Lemma \ref{lm:row_compare} is roughly $\calO(r)$ while the sufficient number of sampled columns in Lemma \ref{lm:columnspmle} is of order $\calO \left( r \frac{N_2}{\underset{i}{\min}\: n_i} \right)$. Thus, the sufficient number of randomly sampled rows is invariant to the distribution of the columns, while the number of columns grows linearly in $N_2/\min_i n_i$, 
and in turn increases if the distribution of the columns is less uniform.}

In Algorithm 3, the columns can admit a clustering structure. Thus, we start the sampling with row sampling.
Let $\hat{r}$ denote a known upper bound on $r$. Such knowledge is often available as side information depending on the particular application. For instance, facial images under varying illumination and facial expressions are known to lie on a special
low-dimensional subspace \cite{wright2009robust}.
For visualization, Fig. \ref{fig: visula 1} provides a simplified illustration of the matrices defined in this section. We sample $C_r \hat{r}$ rows of $\bD$ uniformly at random. Let $\bD_w \in \mathbb{R}^{ (C_r \hat{r}) \times N_2 }$ denote the matrix of sampled rows. We choose $C_r$ sufficiently large to ensure that the non-sparse component of $\bD_w$ is a LR matrix. Define $\bL_w$, assumably with rank $r$, as the LR component of $\bD_w$. 
If we locate a subset of the columns of $\bL_w$ that span its CS, the corresponding columns of $\bL$ would span its CS. To this end, the convex program (\ref{eq2}) is applied to $\bD_w$ to extract its LR component denoted $\hat{\bL}_w$. Then, Algorithm 2 is applied to $\hat{\bL}_w$ to find a set of informative columns 
by sampling $m_1 \approx C_r \hat{r}$ columns.  The matrix $\bD_{s1}$ is formed using the
columns of $\bD$ corresponding to the sampled columns of $\hat{\bL}_w$.
\\
\subsubsection{CS learning}
 \textcolor{black}{Similar to the CS learning step of Algorithm 1, we  obtain the CS of $\bL$ by decomposing $\bD_{s1}$.
Adaptive sampling makes the number of sampled columns required \emph{to capture the CS} of $\bL$ invariant to the distribution of the columns since only informative columns are sampled. Yet, another important advantage of adaptive sampling is that it also makes the sufficient number of sampled columns \emph{to ensure correct decomposition} of $\bD_{s1}$ almost invariant to the distribution of the columns.}\\
\textcolor{black}{To show this fact, we present the following example along with its theoretical footing.
Suppose $\bL$ follows Assumption 1. Thus, per Lemma \ref{lm:row_compare} a small number of randomly sampled rows can capture the RS. Accordingly, the RS  of $\bL_w$ (the LR component of $\bD_w$) is equal to the RS of $\bL$. Thus, if $\bD_w$ is decomposed correctly and Algorithm 2 is applied to $\hat{\bL}_w$, it samples $C_r \frac{r}{n}$ columns from each matrix $\bL_i = \bU_i \bQ_i$ whp. Since Algorithm 2 is deterministic, it is hard to analyze the decomposition algorithm with deterministic sampling and upper bound the coherency parameters. Instead, we consider an imaginary randomized sampling algorithm whose performance emulates that of the proposed sampling approach when $\bL$ follows Assumption 1.}

\smallbreak
\noindent
\textbf{Imaginary sampler:} \textcolor{black}{{Suppose matrix $\bG$ can be represented as $\bG = [\bG_1 \: ... \: \bG_n]$. The CS of matrices $\{ \bG_i \}_{i = 1}^n$ are independent subspaces with dimensions $\{ r_i \}_{i  =1}^n$, respectively. The imaginary sampler applied to $\bG$ samples $C \sum_{i = 1}^n r_i $ consisting of $C r_i$ columns sampled uniformly at random from each submatrix $\{ \bG_i \}_{i=1}^n$.}}

\smallbreak
\noindent

\begin{remark}
The obtained low rank component of $\bD_w$, $\hat{\bL}_w$, is only used to locate the indices of the informative columns. Hence, it is not imperative to obtain exact decomposition of $\bD_w$ and we can still locate the informative columns if $\hat{\bL}_w$ approximates $\bL_w$. However, for analysis purposes in the following lemma it is assumed that $\hat{\bL}_w = {\bL}_w$.

\textcolor{black}{
If the rank of $\bL_w$ is equal to $r$ and $\hat{\bL}_w$ approximates $\bL_w$ well enough, the only difference between the actual and the imaginary sampler is that the former samples the informative columns of $\bL_i$, and the latter samples the columns corresponding to the $i$-th cluster uniformly at random from the columns of $\bL_i$. 
However, sampling random columns from a given cluster is no different that sampling informative columns since $r/n$ randomly sampled columns of $\bL_i$ span its CS whp.} 
\end{remark}

\noindent
\textcolor{black}{The following two lemmas confirm that with adaptive column sampling, the sufficient number of sampled columns to ensure correct decomposition of $\bD_{s1}$ is invariant to the distribution of the columns.
\begin{lemma}
Suppose $\bL$ follows Assumption 1, the imaginary sampler is applied to ${\bL_w}$ to form $\bD_{s1}$, $\bL_{s1}$ has the same CS of $\bL$, and the support set of $\bS$ follows the Bernoulli model with parameter $\rho$. If the rank of $\bL$ is equal to $r$,
\begin{eqnarray}
\begin{aligned}
m_1 \ge    {\frac{r}{\rho}_r} \mu^{''} (\log N_1)^2    \: ,
\end{aligned}
\label{eq:sufflm10}
\end{eqnarray}
and $\rho \leq \rho_s$,
then (\ref{eq20}) yields the exact decomposition
 with probability at least
$ 1    - 2 n (m_1/n)^{-3} - c_8 N_1^{-3} \: , $
where $\varphi_3 = \frac{ \max(r/n , \log m_1/n) }{r/n}$ and
\begin{eqnarray}
\mu^{''} \hspace{-1mm}=\hspace{-.5mm} \max \hspace{-.5mm}\left( \hspace{-.75mm}\frac{c_7 \max (r,\log N_1)}{r} , c_7 \varphi_3 , c_7 \varphi_3 ( \frac{c_9 }{\sqrt{6}} \log N_1)^2\hspace{-.75mm} \right).
\label{eq:muzegdef}
\end{eqnarray}
\label{lm:compare}
\end{lemma}
\begin{lemma}
Suppose $\bL$ follows Assumption 1, the columns of $\bD$ are sampled uniformly at random to form $\bD_{s1}$,   $\bL_{s1}$ has the same CS of $\bL$, and the support set of $\bS$ follows the Bernoulli model with parameter $\rho$.
If the rank of $\bL$ is equal to $r$,
\begin{eqnarray}
\begin{aligned}
m_1 & \ge    {\frac{r}{\rho}_r} \dot{\mu} (\log N_1)^2
\end{aligned}
\label{eq:suff22}
\end{eqnarray}
where
\begin{eqnarray}
\begin{aligned}
& \dot{\mu} =  \max \bigg( \frac{c_7 \max (r,\log N_1)}{r} , 6 \mu_v , 
\mu_v \left(c_9  \log N_1 \right)^2 \bigg) \:,
\\
& \mu_v = \frac{c_7}{n}  \varphi_1 \left(    \frac{N_2}{ \underset{k}{\min} \: n_k } \right) \: , \: \varphi_1 = \frac{\max (r/n\: , \:\log  \underset{i}{\min} \: n_i)}{r/n} \:,
\end{aligned}
\label{eq:defslm11}
\end{eqnarray}
then (\ref{eq20}) decomposes $\bD_{s1}$ correctly with probability at least $1 - c_8 N_1^{-3}  - 2  \sum_{i=1}^n n_i^{-3} $.
\label{lm:compare2}
\end{lemma}
}

\begin{remark}
The sufficient value for $m_1$ in Lemma \ref{lm:compare} is $\calO (r)$, versus $\calO(r \frac{N_2}{ \underset{i}{\min} \: n_i })$ in Lemma \ref{lm:compare2}.
Hence, through adaptive column sampling not only is the number of sampled columns required to capture the CS of $\bL$ invariant to the distribution of the columns, but also the number of sampled columns required for exact decomposition. 
If we sample from a dataset uniformly at random, the resulting data sketches will have the same structure of the data whp. For example, if most of the columns of $\bL$ are aligned along a given direction, then whp most of the randomly sampled columns of $\bL$ will be aligned along that direction as well.
By contrast, adaptive sampling samples the informative columns regardless of the population of the clusters, thereby balances the distribution of the sampled data. 
Thus, the distribution of adaptively sampled columns is closer to a uniform distribution. For instance, the adaptive sampler applied to the data on the right plot of Fig. \ref{fig:innovationfig} samples an equal number of data points from each cluster even though the number of data points in one cluster is notably greater than the other.
This can also be observed by comparing the RS coherency of matrix $\bL_{s1}$ formed with adaptive versus uniform random sampling. The analysis provided in the proof of Lemma \ref{lm:compare} and Lemma \ref{lm:compare2} shows that the RS coherency of $\bL_{s1}$ following adaptive and uniform sampling is roughly $c_7$ and $c_7 \frac{N_2}{\underset{i}{min} \: n_i}$, respectively. Thus, adaptive sampling significantly  improves the coherency of the matrix of sampled columns.
\end{remark}

\subsubsection{Representation matrix learning}
\textcolor{black}{ In this step, unlike Algorithm 1, the rows are not sampled randomly. Instead, we leverage the information embedded in $\bU$ to select the informative rows. Algorithm 2 is applied to $\bU^T$ to locate $m_2 \approx C r$ rows of $\bU$. Thus, we form the matrix $\bD_{s2}$ from the rows of $\bD$ corresponding to the selected rows of $\bU$.} Then, the representation matrix is learned as the optimal point of (\ref{eq23}).
 Subsequently, the LR matrix can be obtained from the learned CS and the representation matrix.
\textcolor{black}{Since in Algorithm 3 it is assumed that the rows do not follow a clustering structure, the analysis of this step is similar to the analysis of the corresponding step in Algorithm 1.}


We can readily state the following theorem which supports the performance of Algorithm 3. In this theorem, it is assumed that the imaginary sampler is used to sample the columns of $\bD$. As mentioned earlier, the performance of the imaginary sampler is \textcolor{black}{approximately} equivalent to our column sampling procedure if $\hat{\bL}_w = \bL_w$, or $\hat{\bL}_w$ approximates $\bL_w$ well enough.
\begin{theorem}
Suppose $\bL$ follows Assumption 1, the support set of $\bS$ follows the Bernoulli model with parameter $\rho$, the rank of $\bL_w$ is equal to $r$, the imaginary sampler is applied to $\bL_w$ to locate the columns forming $\bD_{s1}$, and applied to $\hat{\bU}$ to locate the rows forming $\bD_{s2}$. If  $m_2$ and $\rho$ satisfy (\ref{eq38}) with $\delta^{'} = \delta/n$,
\begin{eqnarray}
\begin{aligned}
 m_1 \ge \max \Bigg(   n \left( \sqrt{ \frac{r}{n}} + \sqrt{2 \log \frac{2 n}{\delta}}  \right)^2 ,  {\frac{r}{\rho}_r} \mu^{''} (\log N_1)^2   \Bigg), \\
\end{aligned}
\label{eq:newm1}
\end{eqnarray}
where  $\{c_i \}_{i = 1}^9$ are constant numbers,
$\kappa= \frac{\log N_1}{r}$, $\beta$ any real number greater than one,
$\mu^{''}$ is equal to (\ref{eq:muzegdef}), and $\varphi_3 = \frac{ \max(r/n , \log m_1/n) }{r/n}$, then Algorithm 3 yields exact decomposition with probability at least
$ 1   - 4 \delta - n (m_1/n)^{-3} - c_8 N_1^{-3} \: , $
provided that $N_1$ is greater than the RHS of the first inequality of (\ref{eq:newm1}).
\label{theorem22}
\end{theorem}

\textcolor{black}{
\subsection{An alternative approach to the CS Learning Step}
In this section, we present an alternative approach to the CS learning step of Algorithm 3. We utilize the information embedded in matrix $\hat{\bL}_w^s$ (the sampled columns of $\hat{\bL}_w$) to obtain $\bU$. }
In particular, if $\bD_w$ is decomposed correctly, the RS of $\hat{\bL}_w^s$ will be the same as that of $\bL_{s1}$ given that the rank of $\bL_w$ is equal to $r$. Let $\bV_{s1}$ be an orthonormal basis for the RS of $\hat{\bL}_w^s$. Thus, to learn the CS of $\bD_{s1}$ we only need to solve
$
\underset{\hat{\bU}}{\min} \| \bD_{s1} -  \hat{\bU} \bV_{s1}^T \|_1 \:.
$

%

\begin{remark}
The convex algorithm (\ref{eq2}) may not always yield accurate decomposition of $\bD_w$ since structured data may not be sufficiently incoherent, suggesting that the decomposition step can be further improved. Let $\bD_w^s$ be the matrix consisting of the columns of $\bD_w$ corresponding to the columns selected from $\hat{\bL}_w$ to form $\hat{\bL}_w^s$. According to our investigations, an improved $\bV_{s1}$ can be obtained by applying the decomposition algorithm presented in \cite{lamport19} to $\bD_w^s$ and use the RS of $\hat{\bL}_w^s$ as an initial guess for the RS of the non-sparse component of $\bD_w^s$. Since $\bD_w^s$ is low-dimensional \textcolor{black}{(roughly $\calO(r) \times \calO(r)$ matrix)}, this extra step is a low complexity operation.
\end{remark}

\subsection{Column/Row sampling from sparsely corrupted data}
\label{sec:both_col_row}
In Algorithm 3, we assumed that the LR component of $\bD_w$ has rank $r$. However, if the rows are not well-distributed,  
a reasonably sized random subset of the rows may not span the RS of $\bL$. Here, we present a sampling approach which can find the informative columns/rows even when  both the columns and the rows exhibit clustering structures such that a small random subset of the columns/rows of $\bL$ cannot span its CS/RS. The algorithm presented in this section (Algorithm 4) can be independently used as an efficient sampling approach from big data. In this paper, we use Algorithm 4 to form $\bD_w$ if both the columns and rows exhibit clustering structures.   

\begin{figure}[t!]
    \includegraphics[width=0.4\textwidth]{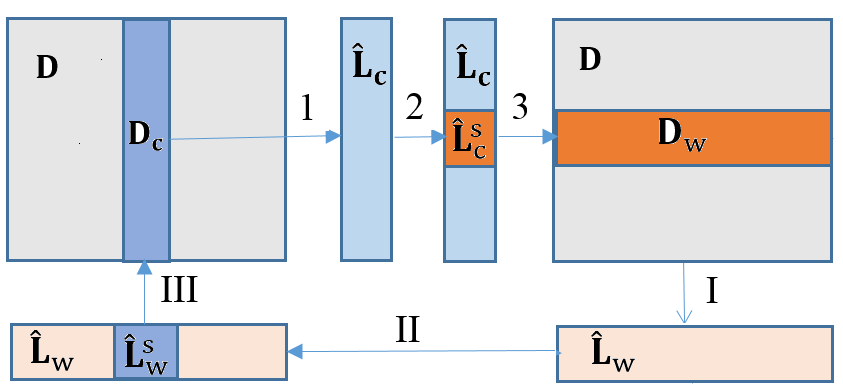}
    \centering
     \vspace{-0.05in}
    \caption{ Visualization of Algorithm 4. We run few cycles of the algorithm and stop when the rank of the LR component does not change over $T$ consecutive steps. One cycle of the algorithm starts from the point marked ``I'' and proceeds as follows. \textbf{I}: Matrix $\bD_w$ is decomposed and $\hat{\bL}_w$ is the obtained LR component of $\bD_w$. \textbf{II}: Algorithm 2 is applied to $\hat{\bL}_w$ to select the informative columns of $\hat{\bL}_w$. $\hat{\bL}_w^s$ is the matrix of columns selected from $\hat{\bL}_w$. \textbf{III}: Matrix $\bD_c$ is formed from the columns of $\bD$ that correspond to the columns of $\hat{\bL}_w^s$.  \textbf{1}: Matrix $\bD_c$ is decomposed and $\hat{\bL}_c$ is the obtained LR component of $\bD_c$. \textbf{2}: Algorithm 2 is applied to $\hat{\bL}_c^T$ to select the informative rows of $\hat{\bL}_c$. $\hat{\bL}_c^s$ is the matrix of rows selected from $\hat{\bL}_c$. \textbf{3}: Matrix $\bD_w$ is formed as the rows of $\bD$ corresponding to the rows used to form $\hat{\bL}_c^s$. }
    \label{alg3}
\end{figure}

The table of Algorithm 4, Fig. \ref{alg3} and its caption provide the details of the proposed sampling approach and the definitions of the used matrices. We start the cycle from the position marked ``I'' in Fig. \ref{alg3} with $\bD_w$ formed according to the initialization step of Algorithm 4. For ease of exposition, assume that $\hat{\bL}_w = \bL_w$ and $\hat{\bL}_c = \bL_c$, i.e., $\bD_w$ and $\bD_c$ are decomposed correctly.
The matrix $\hat{\bL}_w^s$ is the informative columns of $\hat{\bL}_w$. Thus, the rank of $\hat{\bL}_w^s$ is equal to the rank of $\hat{\bL}_w$. Since $\hat{\bL}_w = \bL_w$, $\hat{\bL}_w^s$ is a subset of the rows of $\bL_c$. If the rows of $\bL$ exhibit a clustering structure, it is likely that rank$(\hat{\bL}_w^s)< \text{rank} (\bL_c)$. Thus, rank$({\bL}_w)< \text{rank} (\bL_c)$. 
We continue one cycle of the algorithm by going through steps  1, 2 and 3 of Fig. \ref{alg3} to update $\bD_w$. Using a similar argument, we see that the rank of an updated $\bL_w$ will be greater than the rank of $\bL_c$. Thus, if we run more cycles of the algorithm -- each time updating $\bD_w$ and $\bD_c$ -- the rank of $\bL_w$ and $\bL_c$ will increase.
As detailed in the table of Algorithm 4, we stop if the dimension of the span of the obtained LR component does not change in $T$ consecutive iterations. While there is no guarantee that the rank of $\bL_w$ will converge to $r$ (it can converge to a value smaller than $r$), our investigations have shown that Algorithm 4 performs quite well and the RS of $\bL_w$ converges to the RS of $\bL$ in few steps. We have also found that adding some randomly sampled columns (rows) to $\bD_c (\bD_w)$ can effectively avert converging to a lower dimensional subspace.
For instance, some randomly sampled columns can be added to $\bD_c$, which was obtained by applying Algorithm 2 to $\hat{\bL}_w$.
\begin{algorithm}
\caption{\textcolor{black}{Efficient Column/Row Sampling from Sparsely Corrupted LR Matrices}}
{\footnotesize
\textbf{1. Initialization}\\
Form $\bD_{w} \in \mathbb{R}^{C_r \hat{r} \times N_2}$ by randomly choosing $C_r \hat{r}$ rows of $\bD$. Initialize $k = 1$ and set $T$ equal to an integer greater than 1.
\\
\textbf{2. While} $k > 0$
\\
 \textbf{2.1 Sample the most informative columns}
 \\
\textbf{2.1.1} Obtain $\hat{\bL}_w$  via (\ref{eq2}) as the LR component of $\bD_w$. \\
\textbf{2.1.2} Apply Algorithm 2 to $\hat{\bL}_w$ with $C = C_r$. \\
\textbf{2.1.3} Form the matrix $\bD_c$ from the columns of $\bD$ corresponding to the sampled columns of $\hat{\bL}_w$.
\\
 \textbf{2.2 Sample the most informative rows}
 \\
\textbf{2.2.1} Obtain $\hat{\bL}_c$ via (\ref{eq2}) as the LR component of $\bD_c$. \\
\textbf{2.2.2} Apply Algorithm 2 to $\hat{\bL}_c^T$ with $C = C_r$.\\
\textbf{2.2.3} Form the matrix $\bD_w$ from the rows of $\bD$ corresponding to the sampled rows of $\hat{\bL}_c$.  \\
\textbf{2.3 If} the dimension of the RS of $\hat{\bL}_w$ does not increase in $T$ consecutive iterations, set $k = 0$ to stop the algorithm.
 \\
\textbf{2. End While} \\
\textbf{Output:} The matrices $\bD_w$ and $\hat{\bL}_w$ can be used for column sampling in the first step of the Algorithm presented in Section \ref{sec:alg_eff_sampling}.
}
\end{algorithm}

Algorithm 4 was found to converge in a very small number of iterations (typically less than 4). Thus, even when Algorithm 4 is used
to form the matrix $\bD_w$, the order of complexity of the proposed decomposition method with efficient column/row sampling is $\calO(\max(N_1,N_2) r^2)$.

\section{Online Implementation}
\label{sec:online}
The proposed decomposition approach consists of two main steps, namely, learning the CS of the LR component then decomposing the columns independently. This structure lends itself to online implementation, which could be very beneficial in settings where the data arrives on the fly. The idea is to first learn the CS of the LR component from a small batch of the data and keep tracking the CS. Since the CS is being tracked, any new data column can be decomposed based on the updated subspace. The table of Algorithm 5 details the proposed online matrix decomposition algorithm, where $\bd_t$ denotes the $t^{\text{th}}$ received data column.

Algorithm 5 uses a parameter $n_u$ which determines the rate at which the algorithm updates the CS of the LR component. For instance, if $n_u = 20$, then the CS is updated every 20 new data columns (step 2.2 of Algorithm 5). The parameter $n_u$ has to be set in accordance with the rate of the change of the subspace of the LR component; a small value for $n_u$ is used if the subspace is changing rapidly. The parameter $n_s$ determines the number of columns last received that are used to update the CS. If the subspace changes rapidly, the older columns may be less relevant to the current subspace, hence a small value for $n_s$ is used. On the other hand, when the data is noisy and the subspace changes at a slower rate, a larger value for $n_s$ can lead to more accurate estimation of the CS.

\begin{algorithm}
\caption{Online Implementation}
{\footnotesize
\textbf{1. Initialization} \\
\textbf{1.1} Set the parameters $n_u$ and $n_s$ equal to integers greater than or equal to one.
\\
\textbf{1.2} Form $\bD_0 \in \mathbb{R}^{N_1 \times (C_r \hat{r})}$ as $\bD_0 = [\bd_1 \: \bd_2 \: ... \: \bd_{C_r \hat{r}}].$ Decompose $\bD_0$ using (2) and obtain the CS of its LR component. Define ${\bU}_o$ as the learned CS, ${\bQ}_o$ the appropriate representation matrix and $\hat{\bS}$ the obtained sparse component of $\bD_0$.
\\
\textbf{1.3} Apply Algorithm 2 to ${\bU}_o^T$ to construct the row sampling matrix $\bS_2$.
\\
\\
\textbf{2. For} any new data column $\bd_t$ do\\
\textbf{2.1} Decompose $\bd_t$ as
\begin{eqnarray}
\underset{\hat{\bq}_t}{\min} \| \bS_2^T \bd_{t} - \bS_2^T \bU_o \hat{\bq}_t \|_1 \: ,
\label{eq_q_find}
\end{eqnarray}
and update \\
$\bQ_o \gets [\bQ_o \:  \: \bq_t^{*}]$, $\hat{\bS} \gets [\hat{\bS} \:  \: (\bd_t - \bU_o \bq_t^{*})]$, where $\bq_t^{*}$ is the optimal point of (\ref{eq_q_find}).  \\
\textbf{2.2 If} the remainder of $\frac{t}{n_u}$ is equal to zero, update $\bU_o$ as
\begin{eqnarray}
\underset{\hat{U}_o}{\min} \| \bD_t -  \hat{\bU}_o \bQ_o^t \|_1 \: ,
\label{eq_online_up}
\end{eqnarray}
where $\bQ_o^t$ is the last $n_s \hat{r}$ columns of $\bQ_o$
and $\bD_t$ is the matrix formed from the last $n_s \hat{r}$ received data columns. Apply Algorithm 2 to the new $\bU_o^T$ to update the row sampling matrix $\bS_2$.\\
\textbf{2. End For}\\
\textbf{ Output} The matrix $\hat{\bS}$ as the obtained sparse matrix, $\hat{\bL} = \bD - \hat{\bS}$ as the obtained LR matrix and $\bU_o$ as the current basis for the CS of the LR component.
}
\end{algorithm}

\subsection{Noisy data}
In practice, noisy data can be modeled as
$
\bD = \bL + \bS + \bN \:,
$
where $\bN$ is an additive noise component. 
In \cite{lamport22}, it was shown that the program
\begin{eqnarray}
\begin{aligned}
 \underset{\hat{\bL},\hat{\bS}}{\min} \: \: \lambda\|\hat{\bS}\|_1  + \|\hat{\bL} \|_* \quad \text{s. t.} \quad \big\| \hat{\bL} + \hat{\bS} - \bD \big\|_F \leq \epsilon_n \:, \\
\end{aligned}
\label{eq2noisy}
\end{eqnarray}
can recover the LR and sparse components with an error bound that is proportional to the noise level. The parameter $\epsilon_n$ has to be chosen based on the noise level. This modified version can be used in the proposed algorithms to account for the noise. Similarly, to account for the noise in the representation learning problem (\ref{eq23}), the $\ell_1$-norm minimization problem can be modified as follows:
\begin{eqnarray}
\underset{\hat{\bQ} , \hat{\bE}}{\min} \| \bS_2^T \bD- \bS_2^T \bU \hat{\bQ} - \hat{\bE} \|_1 \quad \text{subject to} \quad  \| \hat{\bE} \|_F \leq \delta_n.
\end{eqnarray}
$\hat{\bE} \in \mathbb{R}^{m_2 \times N_2}$ is used to cancel out the effect of the noise and the parameter $\delta_n$ is chosen based on the noise level \cite{candesstab}.

\section{Numerical Simulations}
\label{sec:sim}
In this section, we present some numerical simulations to study the performance of the proposed randomized decomposition method. First, we present a set of simulations confirming our analysis which established that the sufficient number of sampled columns/rows is linear in $r$. Then, we compare the proposed approach to the state-of-the-art randomized algorithm \cite {new1} and demonstrate that the proposed sampling strategy can lead to notable improvement in performance. We then provide an illustrative example to showcase the effectiveness of our approach on real video frames for background subtraction and activity detection. 
Given the structure of the proposed approach, it is shown that side information can be leveraged to further simplify the decomposition task. In addition, a numerical example is provided to examine the performance of Algorithm 4. Finally, we investigate the performance of the online algorithm and show that the proposed online method can successfully track the underlying subspace.

In all simulations, the Augmented Lagrange multiplier (ALM) algorithm \cite{lamport23,lamport2} is used to solve the optimization problem (\ref{eq2}). In addition, the $\ell_1$-magic routine \cite{lamport57} is used to solve the $\ell_1$-norm minimization problems. It is important to note that in all the provided simulations (except in Section \ref{sec:alternating_alg}), the convex program (\ref{eq2}) that operates on the entire data can yield correct decomposition with respect to the considered criteria. Thus, if the randomized methods cannot yield correct decomposition, it is because they fall short of acquiring the essential information through sampling.

\subsection{Phase transition plots}
In this section, we investigate the required number of randomly sampled columns/rows.
The LR matrix is generated as a product $\bL= \bU_r \bQ_r$, where $\bU_r \in \mathbb{R}^{N_1 \times r}$ and $\bQ_r \in \mathbb{R}^{r \times N_2}$. The elements of $\bU_r$ and $\bQ_r$ are sampled independently from a standard normal $\mathcal{N}(0,1)$ distribution.  The sparse matrix $\bS$ follows the Bernoulli model with $\rho = 0.02$. In this experiment, Algorithm 1 is used and the column/rows are sampled uniformly at random.

\begin{figure}[t!]
    \includegraphics[width=0.50\textwidth]{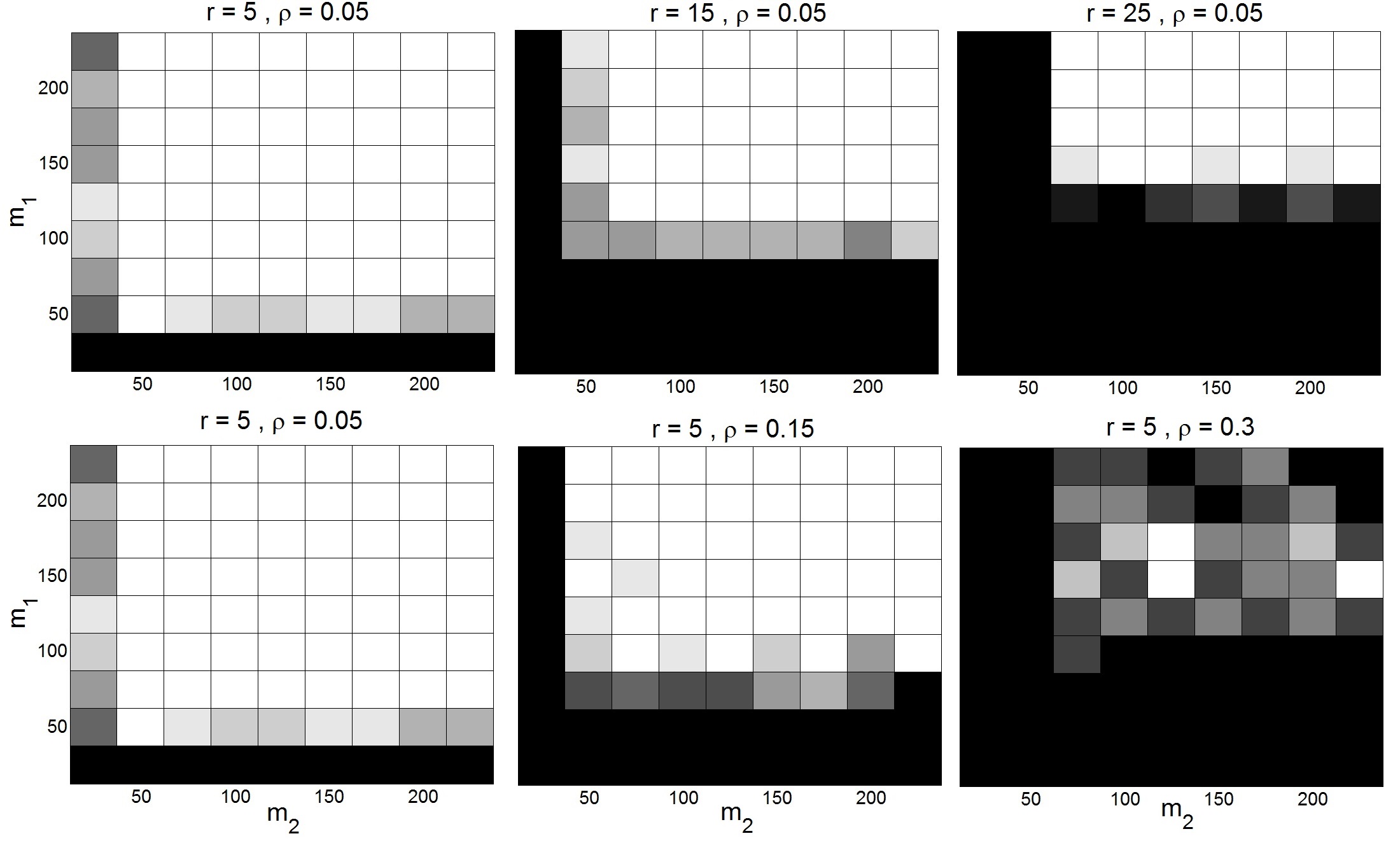}
    \centering
     \vspace{-0.3in}
    \caption{Phase transition plots for various rank and sparsity levels. White designates successful decomposition and black designates incorrect decomposition.}
    \label{fig:trans_plots}
\end{figure}

Fig. \ref{fig:trans_plots} shows the phase transition plots for different numbers of randomly sampled rows/columns. In this simulation, the data is a $1000\times1000$ matrix. For each $(m_1 , m_2)$, we generate 10 random realizations. A trial is considered successful if the recovered LR matrix $\hat{\bL}$ satisfies $\frac{\| \bL - \hat{\bL} \|_F}{\| \bL \|_F} \leq 5 \times 10^{-3}$. It is clear that the required number of sampled columns/rows increases as the rank or the sparsity parameter $\rho$ are increased. When the sparsity parameter is increased to 0.3, the proposed algorithm can hardly yield  correct decomposition. Actually, in this case the matrix $\bS$ is no longer a sparse matrix.

\begin{figure}[t!]
    \includegraphics[width=0.5\textwidth]{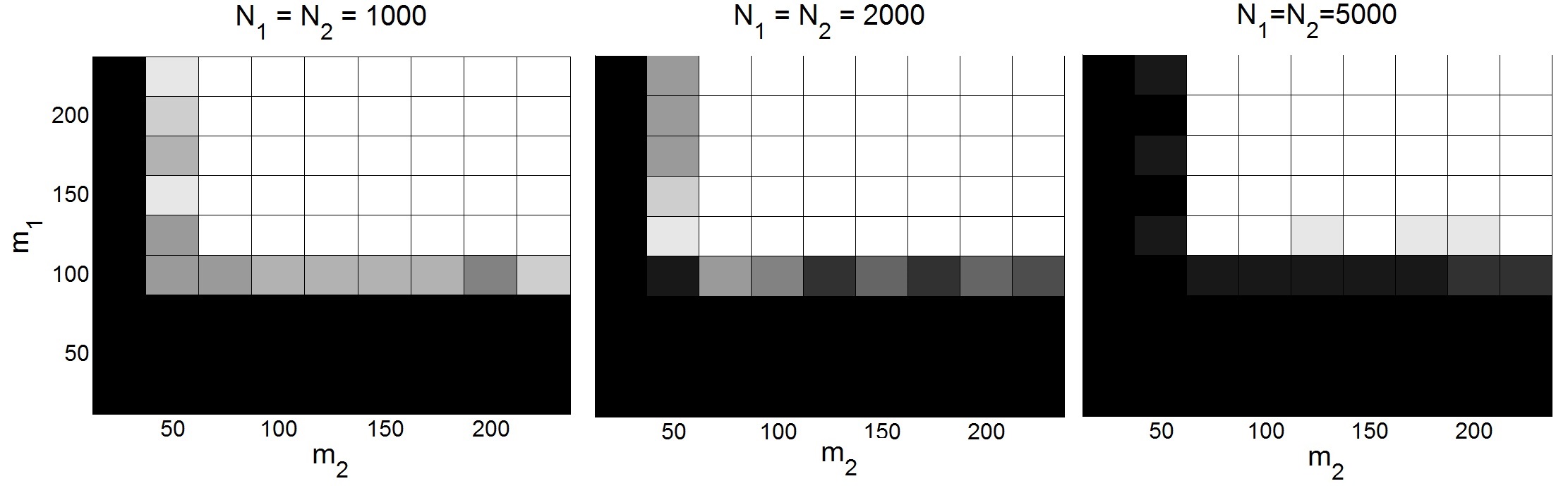}
    \centering
     \vspace{-0.2in}
    \caption{Phase transition plots for various data matrix dimensions ($r = 15 \: , \: \rho = 0.05$).}
    \label{fig:data_dim}
\end{figure}

The top row of Fig. \ref{fig:trans_plots} confirms that the sufficient values for $m_1$ and $m_2$ are roughly linear in $r$. For instance, when the rank is increased from 5 to 25, the required value for $m_1$ increases from 30 to 140. In this experiment, the column and RS of $\bL$ are sampled from the random orthogonal model. Thus, the CS and RS have small coherency whp \cite{lamport4}. Therefore, the important factor governing the sample complexity is the rank of $\bL$. Indeed, Fig. \ref{fig:data_dim} shows the phase transition for different sizes of the data matrix when the rank of $\bL$ is fixed. One can see that the required values for $m_1$ and $m_2$ are almost independent of the size of the data confirming our analysis.

\subsection{Efficient column/row sampling}
In this experiment, Algorithm 3 is compared to the randomized decomposition algorithm in \cite{new1}. It is shown that the proposed sampling strategy can effectively reduce the required number of sampled columns/rows, and makes the proposed method remarkably robust to structured data. In this experiment,  $\bD$ is a $2000\times4200$ matrix. The LR component is generated as $\bL = [ \bG_1 \: \bG_2 \: ... \: \bG_n ] \: . $
For $1 \leq i \leq \frac{n}{2}$, $\bG_i = \bU_i \bQ_i\: ,$
where $\bU_i \in \mathbb{R}^{2000 \times \frac{r}{n}}$, $\bQ_i \in \mathbb{R}^{\frac{r}{n} \times \frac{130 r}{n}}$ and the elements of $\bU_i$ and $\bQ_i$ are sampled independently from a normal distribution $\calN(0 , 1)$. For $n/2+1 \leq i \leq n$, $\bG_i = 13 \bU_i \bQ_i\: ,$
where $\bU_i \in \mathbb{R}^{2000 \times \frac{r}{n}}$, $\bQ_i \in \mathbb{R}^{\frac{r}{n} \times \frac{10 r}{n}}$, and the elements of $\bU_i$ and $\bQ_i$ are sampled independently from an $\calN(0 , 1)$ distribution. We set $r$ equal to 60;  thus, the rank of $\bL$ is equal to 60 whp. The sparse matrix $\bS$ follows the Bernoulli model and each element of $\bS$ is non-zero with probability 0.02. In this simulation, we do not use Algorithm 4 to form $\bD_w$. The matrix $\bD_w$ is formed from 300 uniformly sampled \emph{rows} of $\bD$.

\begin{figure}[t!]
	\centering
    \includegraphics[width=0.5\textwidth]{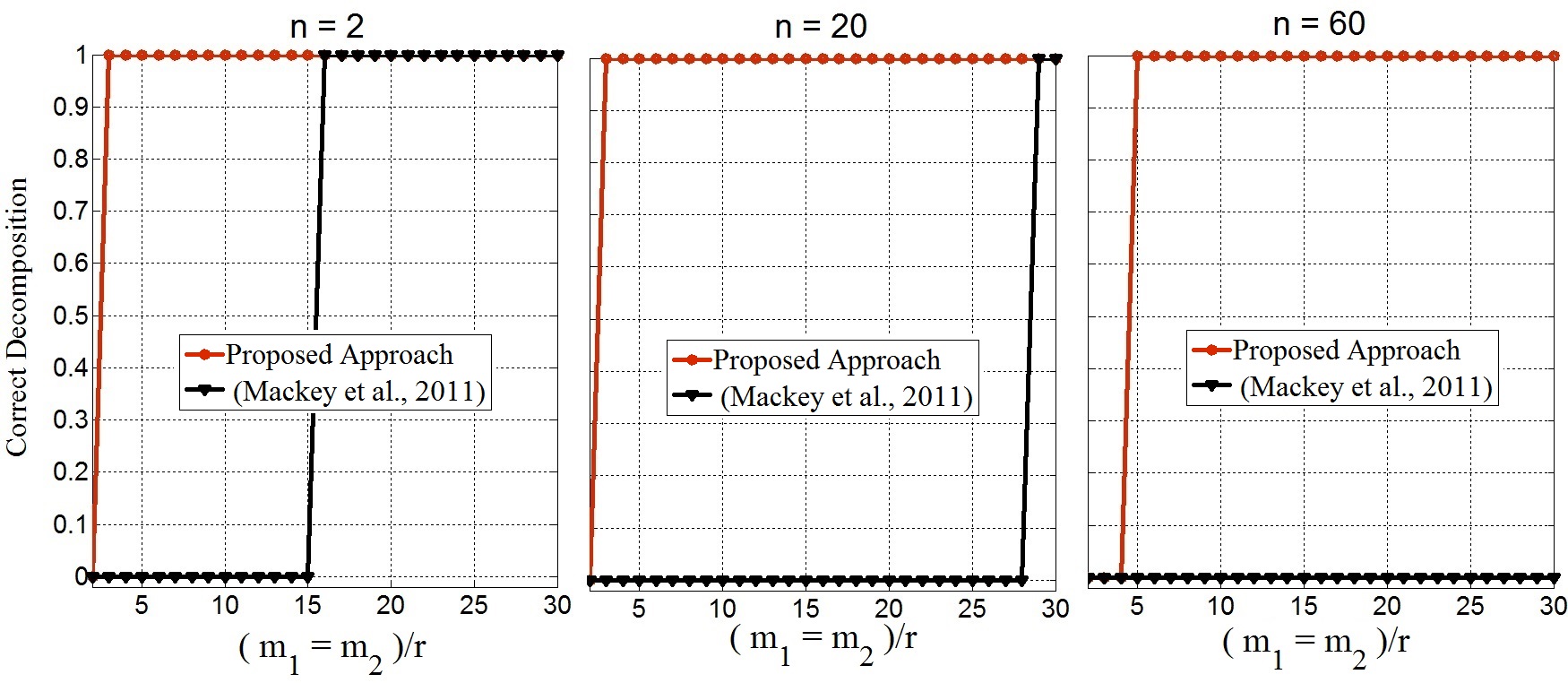}
    \vspace{-.1in}
    \caption{Performance of the proposed approach and the randomized algorithm in \cite{new1}. A value 1 indicates correct decomposition and a value 0 indicates incorrect decomposition.}
    \label{fig:comp}
\end{figure}

 We evaluate the performance of the algorithms for different values of $n$, i.e., different number of clusters. Fig. \ref{fig:comp} shows the performance of the proposed approach and the approach in \cite{new1} for different values of $m_1$ and $m_2$. For each value of $m_1 = m_2$, we compute the error in LR matrix recovery
$\frac{\| \bL - \hat{\bL} \|_F}{\| \bL \|_F} $ averaged over 10 independent runs,
and conclude that the algorithm can yield correct decomposition if the average error is less than 0.01. In Fig. \ref{fig:comp}, the values 0, 1 designate incorrect and correct decomposition, respectively. It can be seen that the presented approach requires a significantly smaller number of samples to yield the correct decomposition. This is due to the fact that the randomized algorithm \cite{new1} samples both the columns and rows uniformly at random and independently. In sharp contrast, we use $\hat{\bL}_w$ to find the most informative columns to form $\bD_{s1}$, and also leverage the information embedded in the CS to find the informative rows to form $\bD_{s2}$. When $n = 60$, \cite{new1} cannot yield correct decomposition even when $m_1 = m_2 = 1800$.

\subsection{Vector decomposition for background subtraction}
\begin{figure}[t!]
	\centering
    \includegraphics[width=0.4\textwidth]{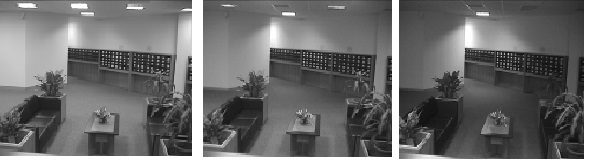}
    \vspace{-.1in}
    \caption{Stationary background.}
    \label{fig:bkgnd}
\end{figure}

The LR plus sparse matrix decomposition can be effectively used to detect a moving object in a stationary background \cite{lamport2,bouwmans2016handbook}. The background is modeled as a LR matrix and the moving object as a sparse matrix. Since videos are typically high dimensional objects, standard algorithms can be quite slow for such applications. Our algorithm is a good candidate for such a problem as it reduces dimensionality significantly. The decomposition problem can be further simplified by leveraging prior information about the stationary background. In particular, we know that the background does not change or we can construct it with some pre-known dictionary. For example, consider the video from \cite{lamport61}, which was also used in \cite{lamport2}. Few frames of the stationary background are illustrated in Fig. \ref{fig:bkgnd}.
Thus, we can simply form the CS of the LR matrix using these frames which can describe the stationary background in different states. Accordingly, we just need to learn the representation matrix. As such, background subtraction is simplified to a vector decomposition problem.
\begin{figure}[t!]
    \includegraphics[width=0.4\textwidth]{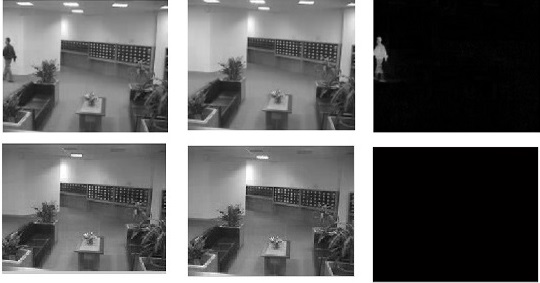}
    \centering
     \vspace{-0.15in}
    \caption{Two frames of a video taken in a lobby. The first column displays the original frames. The second and third columns display the LR and sparse components recovered using the proposed approach. 
 }
    \label{fig:decomp}
\end{figure}
Fig. \ref{fig:decomp} shows that the proposed method successfully separates the background and the moving objects.
In this experiment, 500 randomly sampled rows are used (i.e., 500 randomly sampled pixels) for representation matrix learning (\ref{eq23}). While the running time of our approach is just few milliseconds, it takes almost half an hour if we use (\ref{eq2}) to decompose the video\cite{lamport2}.


\subsection{Alternating algorithm for column sampling}
\label{sec:alternating_alg}
\begin{figure}[t!]
    \includegraphics[width=0.35\textwidth]{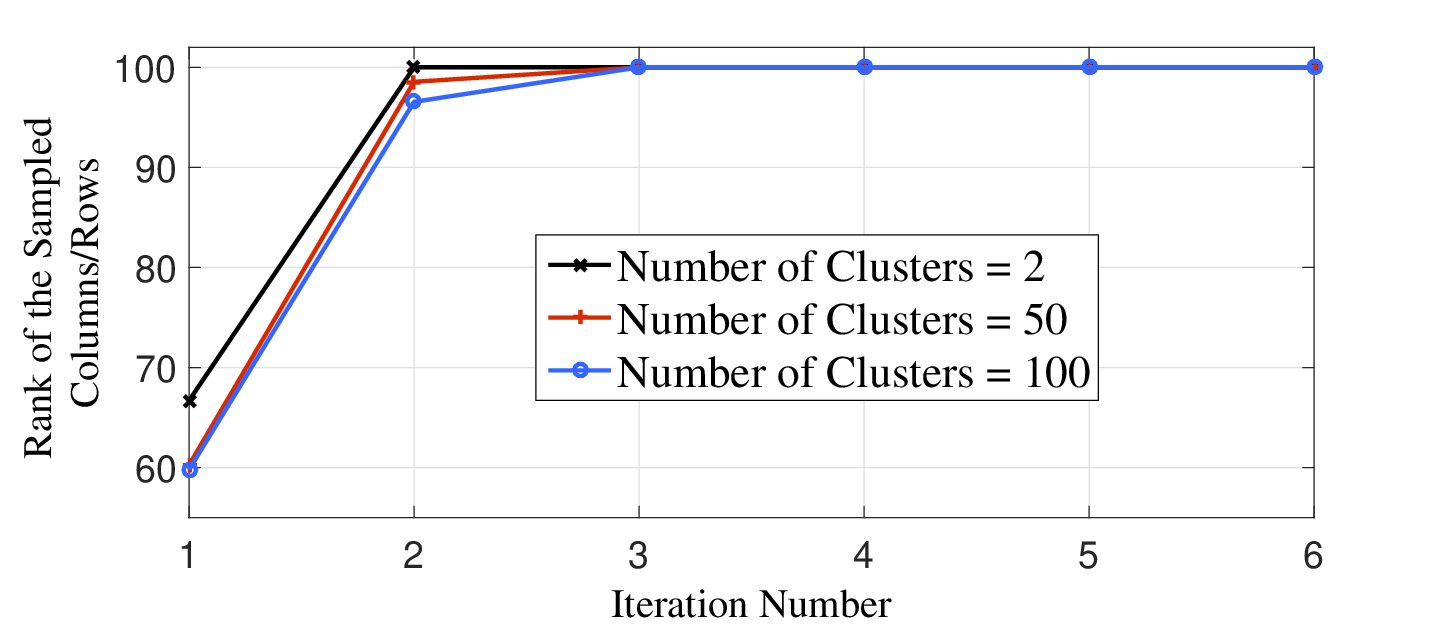}
    \centering
     \vspace{-0.1in}
    \caption{The rank of the matrix of sampled columns.}
    \label{fig:iteration}
\end{figure}

In this section, we investigate the performance of Algorithm 4 for column sampling. The rank of the selected columns is shown to converge to the rank of $\bL$ even when both the rows and columns of $\bL$ exhibit a highly structured distribution.
To generate the LR matrix $\bL$ we first generate a matrix $\bG$ as in Section \ref{sec:non_uniform} but setting $r = 100$. Then, we construct the matrix $\bU_g$ from the first $r$ right singular vectors of $\bG$. We then generate $\bG$ in a similar way and set $\bV_g$ equal to the first $r$ right singular vectors of $\bG$. Let the matrix $\bL = \bU_g \bV_g^T$. For example, for $n=100$, $\bL \in \mathbb{R}^{10250 \times 10250}$. Note that the resulting LR matrix is nearly sparse since in this simulation we consider a very challenging scenario in which both the columns and rows of $\bL$ are highly structured and coherent. Thus, in this simulation we set the sparse matrix equal to zero and use Algorithm 4 as follows. The matrix $\bD_c$ is formed using 300 columns sampled uniformly at random and the following steps are performed iteratively:\\
1. Apply Algorithm 2 to $\bD_c^T$ with $C=3$ to sample approximately $3 r$ columns of $\bD_c^T$ and form $\bD_w$ from the rows of $\bD$ corresponding to the selected rows of $\bD_c$.\\
2. Apply Algorithm 2 to $\bD_w$ with $C=3$ to sample approximately $3 r$ columns of $\bD_w$ and form $\bD_c$ from the columns of $\bD$ corresponding to the selected columns of $\bD_c$. Fig. \ref{fig:iteration} shows the rank of $\bD_c$ after each iteration. It is evident that the algorithm converges to the rank of $\bL$ in less than 3 iterations even for $n = 100$ clusters. For all values of $n$, i.e., $n\in\{2, 50, 60\}$, the data is a $10250\times10250$ matrix.

\subsection{Online Implementation}
\begin{figure}[t!]
    \includegraphics[width=0.5\textwidth]{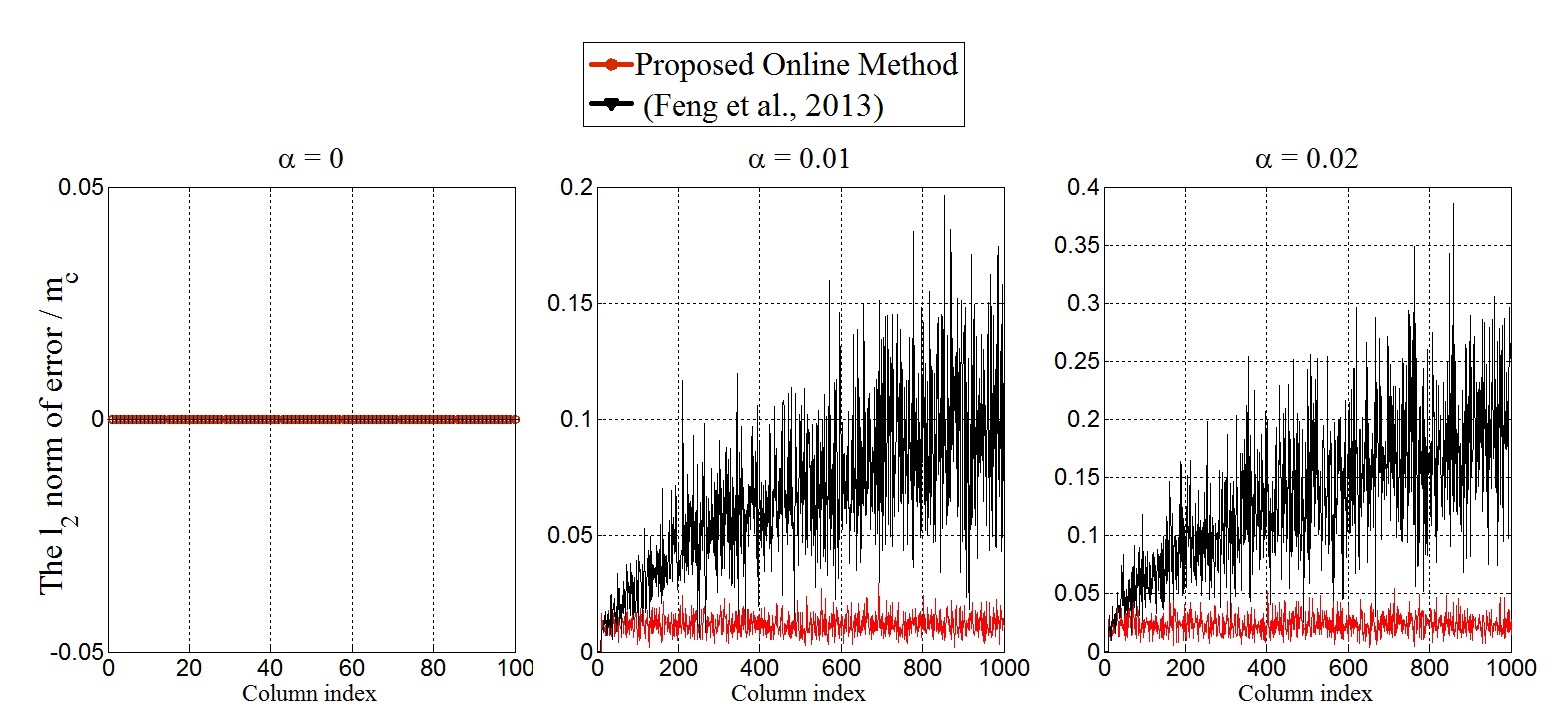}
    \centering
    \vspace{-.1in}
    \caption{Performance of the proposed online approach and the online algorithm in \cite{citonline}. }
    \label{fig: online}
\end{figure}

In this section, the proposed online method is examined. It is shown that the proposed scalable online algorithm tracks the underlying subspace successfully. The matrix $\bS$ follows the Bernoulli model with $\rho = 0.01$. Assume that the orthonormal matrix $\bU \in \mathbb{R}^{N_1 \times r}$ spans a random $r$-dimensional subspace. The matrix $\bL$ is generated as follows.
\smallbreak
\noindent\textbf{For} $k$ from 1 to $N_2$\\
\indent\textbf{1.} Generate $\bE \in \mathbb{R}^{N_1 \times r}$ and  $\bq \in \mathbb{R}^{r \times 1}$  randomly. \\
\indent\textbf{2.} $\bL = [\bL \: \: \bU \bq]  $ . \\
\indent\textbf{3.} \textbf{If} (mod$(k , n)$ = 0) \\
\indent\indent \quad  $  \bU = \text{approx-r} (  \bU + \alpha \bE).$   \\
\indent\indent\textbf{End If} \\
\textbf{End For}
\smallbreak
The elements of $\bq_i $ and $\bE$ are sampled from standard normal distributions. The output of the function approx-r is the matrix of the first $r$ left singular vectors of the input matrix and mod$(k , n)$ is the remainder of $k / n$. The parameters $\alpha$ and $n$ control the rate of change of the underlying subspace. The subspace changes at a higher rate if $\alpha$ is increased or $n$ is decreased. In this simulation, $n = 10$, i.e., the CS is randomly rotated every 10 new data columns. In this simulation, the parameter $r = 5$ and $N_1 = 400$.  We compare the performance of the proposed online approach to the online algorithm in \cite{citonline}. For our proposed method, we set $C=20$ when Algorithm 2 is applied to $\bU$, 
i.e., $20 r$ rows of $\bU$ are sampled. The method presented in \cite{citonline} is initialized with the exact CS and its tuning parameter is set equal to $1/\sqrt{N_1}$. The algorithm \cite{citonline} updates the CS with every new data column. The parameter $n_u$ of the proposed online method is set equal to 4 (i.e., the CS is updated every 4 new data columns) and the parameter $n_s$ is set equal to $5 r$. Define $\hat{\bL}$ as the recovered LR matrix.  Fig. \ref{fig: online} shows the $\ell_2$-norm of the columns of $\bL - \hat{\bL}$ normalized by the average $\ell_2$-norm of the columns of $\bL$ for different values of $\alpha$. One can see that the proposed method can successfully track the CS while it is continuously changing. The online method \cite{citonline} performs well when the subspace is not changing ($\alpha = 0$), however, it fails to track the subspace when it is changing.

\section*{Appendix}
\noindent\textbf{Proof of Lemma \ref{lm3}} \\
The selected columns of $\bL$ can be written as $\bL_{s1}=\bL \bS_1$.
Using the compact SVD of $\bL$, $\bL_{s1}$ can be rewritten as
$\bL_{s1}=\bU \boldsymbol{\Sigma} \bV^T \bS_1$.
Therefore, to show that the CS of $\bL_{s1}$ is equal to that of $\bL$, it suffices to show that the matrix $\bV^T \bS_1$ is a full rank matrix. The matrix $\bS_1$ selects $m_1$ rows of $\bV$ uniformly at random. Therefore, using Theorem 2 in \cite{lamport3}, if
\begin{eqnarray}
m_1 \ge r\gamma^2(\bV) \max \left( c_2 \log r , c_3 \log \frac{3}{\delta} \right),
\label{eq42}
\end{eqnarray}
then the matrix $\bV^T \bS_1$ satisfies the inequality
\begin{eqnarray}
\| I - \frac{N_2}{m_1} \bV^T \bS_1 \bS_1^T \bV \|\leq \frac{1}{2}
\label{eq43}
\end{eqnarray}
with probability at least $(1-\delta)$, where $c_2, c_3$ are numerical constants \cite{lamport3}. Accordingly, if $\sigma_1$ and $\sigma_r$ denote the largest and smallest singular values of $\bS_1^T \bV$, respectively, then
\begin{eqnarray}
\frac{m_1}{2 N_2} \leq \sigma_1^2 \leq \sigma_r^2 \leq \frac{3 m_1}{2 N_2}
\label{eq44}
\end{eqnarray}
Therefore, the singular values of the matrix $\bV^T \bS_1$ are greater than $\sqrt{\frac{m_1}{2 N_2}}$. Accordingly, the matrix $\bV^T \bS_1$ is a full rank matrix.

\begin{remark}
A direct application of Theorem 2 in \cite{lamport3} would in fact lead to the sufficient condition
\begin{eqnarray}
m_1 \ge r\gamma^2(\bR) \max \left( c_2 \log r , c_3 \log \frac{3}{\delta} \right),
\label{neweq42}
\end{eqnarray}
where $\bR \in \mathbb{R}^{N_2 \times N_2}$ denotes the matrix of right singular vectors of $\bL$. The bound in (\ref{eq42}) is slightly tighter since it uses the incoherence parameter $\gamma(\bV)\leq \gamma(\bR) \triangleq \sqrt{N_2} \max_{i,j} |\bR(i,j)|$ in (\ref{neweq42}), where $\bV$ consists of the first $r$ columns of $\bR$. This follows easily by replacing the incoherence parameter in the step that bounds the $\ell_2$-norm of the row vectors of the submatrix in the proof of (\cite{lamport3}, Theorem 2).
%
\end{remark}

\noindent\textbf{Proof of lemma \ref{lm4}}\\
The sampled columns are written as
$
\bD_{s1}=\bD \bS_1= \bL_{s1} + \bS_{s1}.
$
First, we investigate the coherency of the new LR matrix $\bL_{s1}$. Define $\textbf{P}_{\bS_1^T \bV}$ as the projection matrix onto the CS of $\bS_1^T \bV$ which is equal to the rows subspace of $\bL_{s1}$. Therefore, the projection of the standard basis onto the rows subspace of $\bL_{s1}$ can be written as
\begin{eqnarray}
\begin{aligned}
&\underset{i}{\max} \| \textbf{P}_{\bS_1^T \bV} \be_i \|_2^2 = \underset{i}{\max} \| \bS_1^T \bV (\bV^T \bS_1 \bS_1^T \bV)^{-1} \bV^T \bS_1 \be_i \|_2^2 \\
& \leq \underset{j}{\max} \| \bS_1^T \bV (\bV^T \bS_1 \bS_1^T \bV)^{-1} \bV^T  \be_j \|_2^2\\
&\leq \| \bS_1^T \bV (\bV^T \bS_1 \bS_1^T \bV)^{-1} \|^2 \| \bV^T  \be_j \|_2^2 \\
& \leq \frac{\gamma^2 (\bV) r}{N_2} (\frac{\sigma_1^2}{\sigma_r^4})= \frac{\gamma^2 (\bV) r}{N_2} \frac{6 N_2}{m_1}=\frac{(6\gamma^2 (\bV))r}{m_1}
\end{aligned}
\label{eq46}
\end{eqnarray}
where $ (\bS_1^T \bV (\bV^T \bS_1 \bS_1^T \bV)^{-1} \bV^T \bS_1)$ is the projection matrix onto the CS of $\bS_1^T \bV$. The first inequality follows from the fact that $\{ \bS_1 \be_i \}_{i=1}^{m_1}$ is a subset of $\{ \be_j \}_{j=1}^{N_2}$. The second inequality follows from Cauchy-Schwarz inequality and the third inequality follows from (\ref{eq10}) and (\ref{eq44}).

Using lemma 2.2 of \cite{lamport4}, there exists numerical constant $c_7$ such that
$
\underset{i}{\max} \| \bU^T e_i \|_2^2 \leq \frac{\mu_p r}{N_1} \:,
$
with probability at least $1- 2 N_1^{-3}$ and $\mu_p = \frac{c_7 \max (r,\log N_1)}{r}$.\\
In addition, we need to find a bound similar to the third condition of (\ref{eq9}) for the LR matrix $\bL_{s1}$. Let $\bL_{s1} = \bU_{s1} \boldsymbol{\Sigma}_{s1} \bV_{s1}^T$ be the SVD decomposition of $\bL_{s1}$. Define
\begin{eqnarray}
\bH= \bU_{s1} \bV_{s1}^T = \sum_{i=1}^r \bU_{s1}^i (\bV_{s1}^i)^T
\label{eq48}
\end{eqnarray}
where $\bU_{s1}^i$ is the $i^{\text{th}}$ column of $\bU_{s1}$ and $\bV_{s1}^i$ is the $i^{\text{th}}$ column of $\bV_{s1}$.
Given the random orthogonal model of the CS, $\bH$ has the same distribution as
\begin{eqnarray}
\bH^{'} = \sum_{i=1}^r \epsilon_i \bU_{s1}^i (\bV_{s1}^i)^T
\label{eq49}
\end{eqnarray}
where $\{ \epsilon_i \}$ is an independent Rademacher sequence. Using Hoeffding's inequality \cite{lamport43}, conditioned on $\bU_{s1}$ and $\bV_{s1}$ we have
\begin{eqnarray}
\begin{aligned}
& \mathbb{P} \left( |\bH^{'} (i,j)| > t \right) \leq 2 e^{\frac{- t^2}{2 h_{ij}^2}} , \\
& h_{ij}^2=\sum_{k=1}^r (\bU_{s1} (i , k))^2 (\bV_{s1} (j , k))^2.
\label{eq50}
\end{aligned}
\end{eqnarray}
Consider the following lemma adapted from Lemma 2.2 of \cite{lamport4}.
\begin{lemma}[\text{Adapted from lemma 2.2 of \cite{lamport4}}]
If the orthonormal matrix $\bU$ follows the random orthogonal model, then
$
\mathbb{P} \left( |\bU (i,j)|^2 \ge 20\frac{\log N_1}{N_1} \right) \leq 3 N_1^{-8}  .
$
\label{lm5}
\end{lemma}
Therefore,
\begin{eqnarray}
|\bU_{s1} (i , k)|^2 \leq 20 \frac{ \log N_1}{N_1}
\label{eq51}
\end{eqnarray}
with probability at least $1-3 N_1^{-8}$. Thus, we can bound $h_{ij}^2$ as
\begin{eqnarray}
h_{ij}^2 \leq 20\frac{\log N_1}{N_1} \| \bV_{s1} e_i \|_2^2.
\label{eq52}
\end{eqnarray}
Using (\ref{eq46}), (\ref{eq52}) can be rewritten as
\begin{eqnarray}
h_{ij}^2 \leq 120\frac{\log N_1  \gamma^2 (\bV) r}{N_1 m_1}.
\label{eq53}
\end{eqnarray}
Choose $t = \omega \frac{\gamma (\bV) \sqrt{r}}{\sqrt{N_1 m_1}}$ for some constant $\omega$. Thus, the unconditional form of (\ref{eq50}) can be written as
\begin{eqnarray}
\begin{aligned}
&\mathbb{P} \left( |\bH^{'} (i,j)| > \omega \frac{\gamma (\bV) \sqrt{r}}{\sqrt{N_1 m_1} } \right) \leq 2 e^{\frac{- \zeta \omega^2}{\log N_1}} +\\
& \mathbb{P} \left( h_{ij}^2 \ge 120\frac{\log N_1  \gamma^2 (\bV) r}{N_1 m_1} \right)
\end{aligned}
\label{eq54}
\end{eqnarray}
for some numerical constant $\zeta$. Setting $\omega=\zeta^{'} \log N_1$ where $\zeta^{'}$ is a sufficiently large numerical constant gives
\begin{eqnarray}
\mathbb{P} \left( \| \bH^{'} \|_{\infty} \ge  c_9 \log N_1 \frac{\gamma (\bV) \sqrt{r}}{\sqrt{N_1 m_1}} \right) \leq 3 r N_1^{-7}
\label{eq55}
\end{eqnarray}
for some constant number $c_9$ since (\ref{eq51}) should be satisfied for  $r N_1$ random variables. \par
Therefore, according to Lemma \ref{lm1}, if (\ref{eq28}) is satisfied,  the convex algorithm (\ref{eq20}) yields the exact decomposition with probability at least $1-c_8 N_1^{-3}$.
\\
\textbf{Proof of lemma \ref{lm8}}\\
Based on (\ref{eq1}), the matrix of sampled rows can be written as
\begin{eqnarray}
\begin{aligned}
\bD_{s2}=\bS_2^T \bD=\bS_2^T \bL + \bS_2^T \bS
=\bL_{s2} + \bS_{s2}
\end{aligned}
\label{eq62}
\end{eqnarray}
Let $\bL_{s2}=\bU_{s2} \boldsymbol{\Sigma}_{s2} \bV_{s2}^T$ be the compact SVD decomposition of $\bL_{s2}$ and $\bL_{s2}=\bU_{s2}^c \boldsymbol{\Sigma}_{s2}^c (\bV_{s2}^c)^T$ its complete SVD. It can be shown \cite{lamport17} that (\ref{secondrows}) is equivalent to
\begin{eqnarray}
\begin{aligned}
\underset{\hat{\textbf{z}}_{i}}{\min} \quad \| \hat{\textbf{z}}_{i} \|_1  \quad \text{s.t.} \quad
 (\bU_{s2}^{\perp})^T  \hat{\textbf{z}}_{i}= (\bU_{s2}^{\perp})^T \bS_{s2}^i \:.
\end{aligned}
\label{eq63}
\end{eqnarray}
where $\bS_{s2}^i$ is the $i^{th}$ column of $\bS_{s2}$ and $\bU_{s2}^{\perp}$ is the last $(m_2-r)$ columns of $\bU_{s2}^c$ which are orthogonal to $\bU_{s2}$. In other words, if $\bq_i^{*}$ is the optimal point of (\ref{secondrows}) and $\bz_i^{*} \in \mathbb{R}^{m_2}$ is the optimal point of (\ref{eq63}), then $\bz_i^{*} = \bS_2^T (\bd_i - \bU \bq_i^{*})$. Thus, it is enough to show that the optimal point of (\ref{eq63}) is equal to $\bS_{s2}^i$.

 The columns subspace of $\bU_{s2}$ obeys the random orthogonal model. Thus, $\bU_{s2}^{\perp}$ can be modeled as a random subset of $\bU_{s2}^c$. Based on the result in \cite{lamport3}, if we assume that the sign of the non-zero elements of $\bS_{s2}^i$ are uniformly random, then the optimal point of (\ref{eq63}) is $\bS_{s2}^i$ with probability at least $(1-\delta)$ provided that
\begin{eqnarray}
m_2 -r \ge \max \left( c_4 \| \bS_{s2}^i \|_0 \gamma^2 (\bU_{s2}^c) \log \frac{m_2}{\delta} , c_5 \left( \log \frac{m_2}{\delta} \right)^2 \right)
\label{eq64}
\end{eqnarray}
for some fixed numerical constants $c_4$ and $c_5$. The parameter $\gamma (\bU_{s2}^c) = \sqrt{m_2} \underset{i,j}{\max} |\bU_{s2}^c (i,j)|$ and $\| \bS_{s2}^i \|_0$ is the $l_0$-norm of $\bS_{s2}^i$. 
In this paper, we do not assume that the sign of the non-zero elements of the sparse matrix $\bS$ is random. However, according to Theorem 2.3 of \cite{lamport2} (de-randomization technique) if the locations of the nonzero entries of $\bS$ follow the Bernoulli model with parameter $2\rho$, and the signs of $\bS$ are uniformly random and if (\ref{eq63}) yields the exact solution who, then it is also exact with at least the same probability for the model in which the signs are fixed and the locations follow the Bernoulli model with parameter $\rho$ \cite{lamport2}. Therefore, it suffices to provide the sufficient condition for the exact recovery of a random sign sparse vector with Bernoulli parameter $2\rho$.\par
First, we provide sufficient conditions to guarantee that
\begin{eqnarray}
m_2 -r \ge c_4 \| \bS_{s2}^i \|_0 \gamma^2 (\bU_{s2}^c) \log \frac{m_2}{\delta}
\label{eq65}
\end{eqnarray}
with high probability. Using Lemma \ref{lm5} and the union bound,
$
\underset{i,j}{\max} |\bU_{s2}^c (i,j)|^2 \leq 20 \frac{\log m_2}{m_2}
$
with probability at least $1-3m_2^{-6}$.\\

Now, we find the sufficient number of randomly sampled rows, $ m_2$, to guarantee that (\ref{eq65}) is satisfied with high probability. It is obvious that $m_2 < N_1$. Define $\kappa = \frac{\log N_1}{r}$.  Therefore, it is sufficient to show that
\begin{eqnarray}
\frac{m_2}{\|\bS_{s2}^i \|_0}  \ge r \left( c_6  \kappa  \log \frac{N_1}{\delta} +1 \right)
\label{eq67}
\end{eqnarray}
whp, where $c_6=20 c_4 $. Suppose that
\begin{eqnarray}
\rho \leq \frac{1}{\beta r \left( c_6  \kappa \log \frac{N_1}{\delta} +1 \right)}
\label{eq68}
\end{eqnarray}
where $\beta$ is a real number greater than one. Define $\alpha= r \left( c_6  \kappa \log \frac{N_1}{\delta} +1 \right)$. According to (\ref{eq68}) and the Chernoff Bound for Binomial random variables \cite{mcdiarmid1998concentration}, we have
\begin{eqnarray}
\begin{aligned}
& \mathbb{P} \left(  \|\bS_{s2}^i \|_0 -  \frac{m_2}{\beta \alpha} > a \right)
 \leq  \exp \left( \frac{- a^2}{2(\frac{m_2}{\alpha\beta} + \frac{a}{3})}  \right) .
\end{aligned}
\label{eq1n}
\end{eqnarray}
If we set $a=\frac{m_2}{\alpha} \left( 1- \frac{1}{\beta} \right)$, then the inequality (\ref{eq67}) is satisfied. Therefore, (\ref{eq1n}) can be rewritten as
\begin{eqnarray}
\begin{aligned}
& \mathbb{P} \left(  \|\bS_{s2}^i \|_0 -  \frac{m_2}{\beta \alpha} > \frac{m_2}{\alpha} \left( 1- \frac{1}{\beta} \right) \right) \\
&  \leq 2 \exp \left(  \frac{- m_2^2 (\beta -1)^2}{\alpha^2 \beta^2} \frac{3\alpha\beta}{2 m_2(\beta+2)} \right).
\end{aligned}
\label{eq2n}
\end{eqnarray}
Therefore, if
\begin{eqnarray}
m_2 \ge  \frac{ 2 r \beta (\beta - 2) \log \left( \frac{1}{\delta} \right)}{3 (\beta - 1)^2} \left( c_6  \kappa \log \frac{N_1}{\delta} +1 \right) ,
\label{eq3n}
\end{eqnarray}
then the inequality (\ref{eq67}) is satisfied with probability at least $(1-\delta)$. Accordingly, if (\ref{eq38}) is satisfied, then (\ref{eq63}) returns the exact sparse vector with probability at least $1-3 \delta$. The factor 0.5 in the numerator of the RHS of the first inequality of (\ref{eq38}) is due to the de-randomization technique \cite{lamport2} to provide the guarantee for the fixed sign case.
\smallbreak
\noindent\textbf{Proof of Theorem \ref{thm:main_result}} \\
%
The proposed decomposition algorithm yields the exact decomposition if:\\
1. The sampled columns of the LR matrix span the the columns subspace of $\bL$. Lemma \ref{lm3} provides the sufficient conditions on $m_1$ to guarantee that the columns of $\bL_{s1}$ span $\calU$ with high probability. \\
2. The program (\ref{eq20}) yields the correct LR and sparse components of $\bD_{s1}$. Lemma \ref{lm4} provides the sufficient conditions on $m_1$ and $\rho$ to guarantee that $\bD_{s1}$ is decomposed correctly with high probability. \\
3. The sampled rows of the LR matrix span the rows subspace of $\bL$. Since it is assumed that the CS of $\bL$ is sampled from the random orthogonal model, according to Lemma \ref{lm5}
\begin{eqnarray}
\mathbb{P} \left(\max_{i,j} |\bU (i,j)|^2 \ge 20\frac{\log N_1}{N_1} \right) \leq 3 r N_1^{-7} \:.
\label{eq32}
\end{eqnarray}
Therefore, according to Lemma \ref{lm3} if
\begin{eqnarray}
m_2 \ge r \log N_1 \max \left( c_2^{'} \log r , c_3^{'} \log \left(\frac{3}{\delta}\right) \right),
\label{eq33}
\end{eqnarray}
then the selected rows of the matrix $\bL$ span the rows subspace of $\bL$ with probability at least $(1-\delta  - 3 r N_1^{-7})$ where $c_2^{'}$ and $c_3^{'}$ are numerical constants.
\\
4. The minimization (\ref{eq23}) yields the correct RS. Lemma \ref{lm8} provides the sufficient conditions to ensure that (\ref{secondrows}) yields the correct representation vector. In order to guarantee the performance of (\ref{eq23}), we substitute $\delta$ with $\delta/N_2$ since (\ref{eq23}) has to return exact representation for the columns of $\bD$.
Therefore,
\[
\mathbb{P} \left( \text{Incorrect Decomposition} \right) \leq \delta + c_8 N_1^{-3} + \delta + 3r N_1^{-7} + 3\delta.
\]

\smallbreak
\noindent
\textbf{Proof of Lemma \ref{lm:coherece_study}}\\
Since the rank of $\bL$ is equal to $r$, the column spaces of $\{ \bU_i \}_{i = 1}^n$ are independent $r/n$-dimensional subspaces and the rank of the matrices $\{ \bQ_i \}_{i = 1}^n$ is equal to $r/n$. The matrix $\bL$ can be expressed as  $\bL = [ \bU_1 \: ... \: \bU_n] \bO$, where $\bO$ is a block matrix with blocks equal to the matrices $\{ \bQ_i \}_{i = 1}^n$.
  Suppose $ {\bQ_k^o}^T \in \mathbb{R}^{n_k \times r/n  }  $ is
an orthonormal matrix such that the row space of $\bQ_k^o$ is equal to the row space of $\bQ_k$. Thus, according to Lemma \ref{lm5},
\begin{eqnarray}
\mathbb{P} \left( |\bQ_k^o (i,j)|^2 > 20\frac{\log n_k}{n_k} \right) \leq 3 n_k^{-8}  .
\end{eqnarray}
Thus, $\max_{i,j} |\bQ_k^o(i,j)|^2 \leq 20\frac{\log n_k}{n_k}$ with probability at least $1 -  \frac{3r}{n} n_k^{-7}$. Given the block diagonal structure of $\bO$,
\begin{eqnarray}
\begin{aligned}
 \gamma^2(\bV) &= N_2 \max_{k} \left[ \max_{i,j} |\bQ_k^o(i,j)|^2  \right] \\
 & \leq 20 \log (\underset{i}{\max}\: n_i) \frac{ N_2}{ \underset{i}{\min}\: n_i }
\end{aligned}
\label{eq:cohtotal12}
\end{eqnarray}
with probability at least $1 - \frac{3r}{n} \sum_{i=1}^n n_i^{-7}$.
Now we compute $\underset{i}{\max} \| \bV^T \be_i \|_2^2$. According to Lemma 2.2 in \cite{lamport4},
\begin{eqnarray}
 \mathbb{P}\left[ \underset{i}{\max} \: \| \bQ_k^o \be_i \|_2^2 > \frac{c_7 \max(r/n , \log n_k) }{n_k} \right] < 2 n_k^{-3} \:,
\label{eq:coh_ran2}
\end{eqnarray}
where $c_7$ is a numerical constant. Thus, based on (\ref{eq:coh_ran2}) and the block structure of $\bO$,
\begin{eqnarray}
\begin{aligned}
 & \mathbb{P} \Bigg[ \underset{i}{\max} \: \|\be_i^T \bV \|_2^2 =
\underset{k}{\max} \left( \underset{i}{\max} \: \| \bQ_k^o \be_i \|_2^2 \right) > \\
& \frac{c_7 r \varphi_1 }{ N_2 } \left(  \frac{1}{n}  \frac{N_2}{ \underset{k}{\min} \: n_k } \right)  \Bigg] \leq 2  \sum_{i=1}^n n_i^{-3}
\end{aligned}
\label{eq:totalcoh2}
\end{eqnarray}
where $\varphi_1 = \frac{\max (r/n , \log \underset{k}{\max} \: n_k)}{r/n}.$

\smallbreak
\noindent
\textbf{Proof of Lemma \ref{Lm:lowerbound}}\\
The proof of this lemma is similar to the proof of Lemma \ref{lm:coherece_study}. We make use of the following lemma to establish a lower bound on the coherency  of the block diagonal matrix $\bO$.
The proof of Lemma \ref{lm:loweasli} is provided at the end of the Appendix.
\begin{lemma}
Suppose $\bU \in \mathbb{R}^{N_1  \times r}$ is an orthonormal matrix which spans a random $r$-dimensional matrix. If $r \ge 18 \log N_1$ and $N_1 \ge 96 r \log N_1$, then
$
 \mathbb{P} \left[ \underset{i}{\max} \: \| \bU^T \be_i \|_2^2  <  \frac{0.5 \: r}{N_1} \right] \leq 2 N_1^{-5} \: .
$

\label{lm:loweasli}
\end{lemma}

\smallbreak
\noindent
\textbf{Proof of Lemma  \ref{lm:columnspmle}}\\
To prove Lemma \ref{lm:columnspmle}, it suffices to ensure that the number of sampled columns corresponding to each submatrix $\{ \bL_i = \bU_i \bQ_i \}_{i  = 1}^n$ is sufficiently large to span the column space of $\{\bL_i\}_{i  = 1}^n$.
Suppose $ {\bQ_k^o}^T \in \mathbb{R}^{n_k \times r/n  }  $ is
an orthonormal matrix and the RS of $\bQ_k^o$ is equal to the RS of $\bQ_k$.
According to Lemma 2.2 in \cite{lamport4},
\begin{eqnarray}
 \mathbb{P}\left[ \underset{i}{\max} \: \|\bQ_k^o \be_i \|_2^2 > \frac{c_7 \max(r/n , \log n_k) }{n_k} \right] < 2 n_k^{-3} \: .
\label{eq:coh_ran22}
\end{eqnarray}
Define $\bL_i^s$ as the columns sampled from submatrix $\bL_i$. According to Lemma 5 in \cite{new3} and (\ref{eq:coh_ran22}), if the number of columns of $\bL_i^s$ is greater than or equal to
$
\xi_i = {10}\: c_7 \: {\max(r/n , \log n_i)} \log \frac{2 r}{\delta} \:,
$
then the CS of $\bL_i^s$ is equal to the CS of $\bL_i$ with probability at least $1 - \delta/n - 2 n_k^{-3}$.
According to Lemma 6 in \cite{new3}, if
$
m_1 \ge \left( 2 + \frac{3}{\xi_i} \log\frac{2 n}{\delta} \right)  \frac{\xi_i N_2}{ n_i } \:,
$
then the number of columns sampled from $\bL_i$ is greater than $\xi_i$ with probability at least $1- \delta/n$. Thus, if (\ref{eq:suf_sample_c}) is satisfied, the sampled columns span the CS of $\bL$ with probability at least $1 - 2 \delta - 2 \sum_{i = 1}^n n_i^{-3}$.


\smallbreak
\noindent
\textbf{Proof of Lemma \ref{lm:row_compare}}\\
From Lemma 2.2 in \cite{lamport4},
\begin{eqnarray}
 \mathbb{P}\left[ \underset{i}{\max} \: \|\be_i^T \bU \|_2^2 > \frac{c_7 \max(r , \log N_1) }{N_1} \right] < 2 N_1^{-3} \:,
\end{eqnarray}
where $\bU$ is an orthonormal basis for the CS of $\bL$. Thus, according to Lemma 5 in \cite{new3}, if (\ref{eq:suff_row_com}) is true, the sampled rows span the RS of $\bL$ with probability at least $1 - 2 N_1^{-3} - \delta$.

\smallbreak
\noindent
\textbf{Proof of Lemma \ref{lm:compare}} \\
Define $\bQ_i^s$ as the columns of $\bQ_i$ corresponding to the sampled columns from $\bL_i = \bU_i \bQ_i$.
First, we establish an upper bound on the RS coherency of the matrix of sampled columns. Define $\bP_{{Q_i^s}^T}$ as the projection matrix onto the row space of $\bQ_i^s$. From Lemma 2.2 in \cite{lamport4},
\begin{align}
 \mathbb{P} \left[ \underset{i}{\max} \| \bP_{{Q_i^s}^T} \be_i\|_2^2 > \frac{c_7  \max (r/n , \log Cr/n)}{Cr/n} \right] <
 2 \left(  C \frac{ r}{n} \right)^{-3}.
\end{align}
Since the rank of $\bL$ is equal to $r$, the matrices $\{ \bU_i \}_{i = 1}^n$ span $n$ independent $r/n$-dimensional linear subspaces.
Thus, similar to the analysis used in the proof of Lemma \ref{lm:columnspmle},
\begin{eqnarray}
\begin{aligned}
 \mathbb{P}  \left[ \underset{i}{\max} \| \textbf{P}_{\bS_1^T \bV} \be_i \|_2^2 >
 \frac{c_7 r \varphi_3}{m_1} \right] < 2 n \left(  \frac{C r}{n}  \right)^{-3} \:,
 \end{aligned}
\end{eqnarray}
where $\textbf{P}_{\bS_1^T \bV}$ is the projection matrix onto the RS of $\bL_{s1}$, $m_1$ is the number of columns sampled to form $\bD_{s1}$ (which in this lemma is equal to $C r$), and
$\varphi_3 = \frac{ \max(r/n , \log Cr/n) }{r/n}.$

Finally, similar to the analysis provided in the proof of Lemma \ref{lm4},  if  (\ref{eq:sufflm10}) is satisfied,
then the convex algorithm (\ref{eq20}) yields exact decomposition with probability at least $1-2 N_1^{-3} - 2 n  \left( \frac{m_1}{n} \right)^{-3}$ where
$\mu^{''}$ is equal to (\ref{eq:muzegdef}).

\smallbreak
\smallbreak
\noindent
\textbf{Proof of Lemma \ref{lm:compare2}}
\\
According to Lemma \ref{lm:coherece_study},
\begin{eqnarray}
\begin{aligned}
 & \mathbb{P} \Bigg[ \underset{i}{\max} \: \|\be_i^T \bV \|_2^2  > \frac{c_7 r \varphi_1 }{ N_2 } \left(  \frac{1}{n}  \frac{N_2}{ \underset{k}{\min} \: n_k } \right)  \Bigg] \leq 2  \sum_{i=1}^n n_i^{-3} \\
\end{aligned}
\label{eq:chostudy}
\end{eqnarray}

\noindent
Thus, based on (\ref{eq:chostudy}) and the analysis provided in the proof of Lemma \ref{lm4}, if (\ref{eq:suff22}) is satisfied,
then (\ref{eq20}) decomposes $\bD_{s1}$ correctly
 with probability at least $1 - 2 N_1^{-3} - 2  \sum_{i=1}^n n_i^{-3}$.
 
\smallbreak
\noindent
\textbf{\textcolor{black}{Proof of Theorem} \ref{theorem22}}\\
The imaginary sampler samples $m_1 = C r$ columns of $\bD$.
First, we ensure that the sampled columns corresponding to each submatrix $\bL_i = \bU_i \bQ_i$ span the column space of $\bL_i$. Define $\bQ_i^s$ as the columns of $\bQ_i$, corresponding to the columns sampled from $\bL_i$. Based on the following lemma from \cite{vershynin2010introduction,davidson2001local}, if
$
\frac{C r}{n} \ge \left( \sqrt{ \frac{r}{n}} + \sqrt{2 \log \frac{2 n}{\delta}}  \right)^2
$,
then the rank of $\bQ_i^s$ is equal to $\frac{r}{n}$ with probability at least $1 - \delta/n$.
\begin{lemma}
Let $\bA$ be an $N \times n $ matrix whose entries are independent standard normal variables. Then for every $t \ge \sqrt{2 \log 2/\delta} $,
\begin{eqnarray}
\sqrt{N} - \sqrt{n} - t \leq \sigma_{min} (\bA) \leq \sigma_{max} (\bA) \leq \sqrt{N} + \sqrt{n} + t \:
\end{eqnarray}
with probability at least $1 - \delta$, where $\sigma_{min} (\bA)$ and $\sigma_{max} (\bA)$ are the minimum and maximum singular values of $\bA$. 
\label{lm:singulars}
\end{lemma}

\noindent
Lemma \ref{lm:compare} establishes a sufficient condition to guarantee the performance of the CS learning step. In addition, since the rows lie in one subspace and the imaginary sampler samples the rows within each subspace uniformly at random, the analysis of the representation learning step is similar to the analysis of the corresponding step of Algorithm 1.

\smallbreak
\noindent
\textbf{Proof of Lemma \ref{lm:loweasli}}\\
First we review the result below from \cite{laurent2000adaptive}.
\begin{lemma}
Let $Y_r$ be a chi-squared random variable with $r$ degrees of freedom. Then
for each $t > 0$ and for each $\epsilon \in (0 , 1)$
\begin{eqnarray}
\begin{aligned}
& \mathbb{P} \left[ Y_r - r \leq -t \sqrt{2 r} \right] \leq e^{- t^2/2} \\
& \mathbb{P} \left[ Y_r \ge r (1- \epsilon)^{-1} \right] \leq e^{- \epsilon^2 r /4 } .
\end{aligned}
\label{eq:lmchis}
\end{eqnarray}
\label{lm:chisqq}
\end{lemma}
The distribution of $\| \bU^T \be_i \|_2^2$ is equivalent to the distribution of $\frac{Y_r}{Y_{N_1}}$ \cite{lamport4}. For each $\lambda < 0 $, it follows from Lemma \ref{lm:chisqq} that
\begin{eqnarray}
\begin{aligned}
& \mathbb{P} \left[ \| \bU^T \be_i \|_2^2 - \frac{r}{N_1} \leq \frac{\lambda \sqrt{2 r} }{N_1} \right]  = \mathbb{P} \left[ Y_r \leq (r + \lambda \sqrt{2 r}) \frac{Y_{N_1}}{N_1}  \right] \\
& \leq \mathbb{P} \left[ Y_r \leq (r + \lambda \sqrt{2r} )(1 - \epsilon)^{-1} \right] + e^{- \epsilon^2 N_1 /4}  \:,
\end{aligned}
\label{eq:prob_1}
\end{eqnarray}
where the inequality follows from the second inequality of (\ref{eq:lmchis}).
Set $\epsilon= \frac{1}{1 + \sqrt{r}}$ and $\lambda = -3/2 \sqrt{\log N_1}$. If $r > 1$, then $(1- \epsilon)^{-1} \leq 2$. Thus, from Lemma \ref{lm:chisqq}, the first term of the RHS of (\ref{eq:prob_1}) can be expanded as
\begin{eqnarray}
\begin{aligned}
\mathbb{P}  \left[ Y_r - r \leq \sqrt{2r} \left( \frac{1}{\sqrt{2}} - 3 \sqrt{\log N_1}   \right)   \right] 
&\leq e^{ - 6 \log N_1 } = N_1^{-6}.
\end{aligned}
\label{eq:prob_2}
\end{eqnarray}
Accordingly,
$
 \mathbb{P} \left[ \| \bU^T \be_i \|_2^2 - \frac{r}{N_1} \leq \frac{\lambda \sqrt{2 r} }{N_1} \right] 
 \leq N_1^{-6} + e^{\frac{- N_1}{ 4 (\sqrt{r} +1)^2 }} \leq N_1^{-6} + e^{\frac{-N_1}{16 r}}.
$
Thus, if $N_1 \ge 96 r \log N_1$ and $r \ge 18 \log N_1$, then
\begin{eqnarray}
\begin{aligned}
 \underset{i}{\max} \: \| \bU^T \be_i \|_2^2 \ge \frac{r}{N_1} -  \frac{3/2 \sqrt{ 2 r \log N_1 } }{N_1}  \ge \frac{0.5 \: r}{N_1}
\end{aligned}
\end{eqnarray}
with probability at least $1 - 2 N_1^{-5}$.

\bibliographystyle{IEEEtran}
\bibliography{IEEEabrv,bibfile}

\end{document}